\documentclass [12pt] {amsart}
\usepackage{enumitem}
\usepackage{times}
\usepackage{latexsym,amssymb, amsmath}
\usepackage{graphicx}
\usepackage{pifont}

\usepackage{simplemargins}
\setleftmargin{2.75cm}
\setrightmargin{2.75cm}
\settopmargin{2.0in}
\setbottommargin{2.0in}
\setallmargins{2.75cm}

\def\pf{\noindent\emph{Proof: }}
\def\stop{\hfill$\Box$}

\usepackage{amsmath, amsfonts, amssymb, amsthm}
\newtheorem{thm}{Theorem}
\newtheorem{cor}[thm]{Corollary}
\newtheorem{lemma}[thm]{Lemma}

\newtheorem{defn}[thm]{Definition}
\newtheorem{prop}[thm]{Proposition}

\numberwithin{thm}{section}

\DeclareMathOperator{\sgn}{sgn}

\begin{document}

\title{On the Local Regularity of the Hilbert Transform}

\author{Yifei Pan}
\address{Department of Mathematical Sciences\\
	Purdue University Fort Wayne\\
	Fort Wayne, Indiana 46805}
\email{pan@pfw.edu}

\author{Jianfei Wang}
\address{School of Mathematical Sciences\\
	Huaqiao University \\ Quanzhou, Fujian 362021 \\China
	}
\email{wangjf@mail.ustc.edu.cn}

\author{ Yu Yan  }
\address {Department of Mathematics and Computer Science\\ Biola University\\La Mirada, California 90639}
\email {yu.yan@biola.edu}

\newtheorem{Thm}{Theorem}
\newtheorem{Lemm}{Lemma}
\newtheorem{Cor}{Corollary}

\begin{abstract}
In this paper the local regularity of the Hilbert transform is studied, and local smoothness and real analyticity results are obtained.
\end{abstract}

\maketitle

\section{Introduction}

\vspace{.1in}
Given a function $f \in L^p(\mathbb{R}) (1 \leq p < \infty)$, its Hilbert transform, denoted $Hf$, is well defined a.e. as the principal value of the convolution of $f$  with the function $\frac{1}{\pi x} $:

\begin{equation*}
	\label{eqn:Hilbert}
	Hf(x) \,\, = \,\,  p.v. \, \frac{1}{\pi} \int _{-\infty}^{\infty} \frac{f(t)}{x-t} \, dt \,\, := \,\, \lim _{\epsilon \to 0} \frac{1}{\pi} \int _{|x-t| > \epsilon} \frac{f(t)}{x-t} \, dt.  
\end{equation*}

\vspace{.1in}

The purpose of this paper is to study the local regularity of the Hilbert transform when only the local regularity of the function is assumed, and the main result is the following theorem.

\begin{thm}
	\label{thm:main}
	Suppose $f \in L^p(\mathbb{R}) (1 \leq p < \infty)$.  If $f$ is $C^{k, \alpha} (k \geq 0, 0<\alpha<1)$ in a neighborhood of a point $x_0 \in \mathbb{R}$, then its Hilbert transform $Hf$ is $C^{k, \alpha_0}$ in a neighborhood of $x_0$ for any $0<\alpha_0 < \alpha$. 
\end{thm}

\noindent
Here we use the standard notations $C^{0,\alpha}$ for functions that are H\"older continuous with exponent $\alpha$ in a bounded interval and $C^{k,\alpha}$ for functions whose $k$-th derivatives are H\"older continuous with exponent $\alpha$.

\vspace{.1in}

It is well known that derivative formulas exist for globally defined differentiable functions; for example (\cite{Pandey}), if $f \in C^{k+1}(\mathbb{R})$ and all its derivatives are in $L^p(\mathbb{R})$, then its Hilbert transform satisfies
$$ \left ( Hf \right )^{(k)} (x) = H \left ( f^{(k)} \right ) (x).   $$

\vspace{.1in}
\noindent
In the case of local differentiability the following corollary is an immediate consequence of Theorem \ref{thm:main}.

\vspace{.05in}

\begin{cor}
	\label{cor:C_infinity}
	Given a function $f(x) \in L^p(\mathbb{R}) (1 \leq p < \infty)$, if $f$ is $C^{k+1} $ near a point, then its Hilbert transform $Hf(x)$ is $C^{k} $ near the same point. In particular, if $f(x)$ is $C^{\infty}$ at a point, then $Hf(x)$ is also $C^{\infty}$ at the same point.	
\end{cor}

\noindent

\vspace{.1in}
\noindent
In addition to local differentiability, the Hilbert transform also preserves local real analyticity.

\vspace{.05in}

\begin{thm}
	\label{thm:real_analytic}
	If a function $f(x)$ is in $L^p(\mathbb{R}) (1 \leq p < \infty)$ and is real analytic near a point, then its Hilbert transform $Hf(x)$ is also real analytic near the same point.	
\end{thm}

\vspace{.1in}
\noindent

 For functions defined on the unit circle, a classical theorem by Privalov (\cite{Privalov_1}, \cite{Privalov_2}, \cite{Bao-Eben-Roth}) states that the Hilbert transform on $S^1$ is a bounded operator from $C^{0, \alpha}(S^1)$ $(0 < \alpha <1)$ into itself.
 In an attempt to generalize Privalov's theorem to the Hilbert transform on $\mathbb{R}$, Dyachenko, Nursultanov, and Tikhonov (\cite{Dyach_NUrs_Tikh}) considered the global H\"older regularity of the Hilbert transform on  a particular class of functions that includes $ C^{0, \alpha} (\mathbb{R}) \cap L^p(\mathbb{R}) $, and they showed that the Hilbert transform is well defined everywhere for this class, furthermore it has bounded  global $C^{0, \alpha}$ norm.  Our Theorems \ref{thm:main} to \ref{thm:real_analytic}, on the other hand, are completely local in the sense that we only assume the function is regular near a point.




\vspace{.1in}
Theorem \ref{thm:main} guarantees that the (possible) regularity loss under the Hilbert transform must be arbitrarily small if a loss does occur; namely, if the function is $C^{k, \alpha}$ near a point, then $Hf$ is at least $C^{k, \alpha_0}$ for any $0 < \alpha_0 < \alpha$ near that point. As will be shown in Section \ref{section:example}, it is possible for local regularity to remain unchanged under the Hilbert transform. Thus a rather interesting question is whether Theorem \ref{thm:main} can be improved to $\alpha_0=\alpha$, but we have found no success in doing so.

\vspace{.1in}
The tool we shall use to deal with the local differentiability  is the so-called finite Hilbert transform.  Given a function $f$ defined on a finite interval $(a,b)$, its finite Hilbert transform is
\begin{equation*}
	\label{eqn:local_Hilbert}
	p.v. \, \frac{1}{\pi} \int _{a}^{b} \frac{f(t)}{x-t} \, dt, 
\end{equation*}
where $a<x<b$.  This is equivalent to $\chi_{(a,b)}(x)  H \left (\chi_{(a,b)}f \right )(x)$, where $\chi$ is the characteristic function. 
In order to work with this singular integral, we introduce the following (non-singular) operator $T$ for $f \in C^{0,\alpha}([a,b])$:
$$Tf(x) = \int_a ^b \dfrac{f(y)-f(x)}{y-x} \, dy.$$

\vspace{.05in}
\noindent
In Section \ref{section:study_T} we will study some fundamental properties of this operator $T$ pertinent to the proofs of Theorems \ref{thm:main} and \ref{thm:real_analytic}. 


\vspace{.1in}

On a side note, the local integrability of the Hilbert transform was considered in the literature.  As a matter of fact, there is essentially a characterization of local integrability.  For example, according to a result by Calder\'on and Zygmund (\cite{Javier}), if $I_1 \subset I_2$ are two open intervals and $f\log^{+}|f| \in L^1(I_2)$, then $Hf \in L^1(I_1)$.  Conversely, Stein (\cite{Stein}, \cite{Javier}) showed that if $f \geq 0$ and $Hf \in L^1(I_2)$, then $f\log^+|f| \in L^1(I_1)$. 
For reference on recent works where the Hilbert transform is exclusively used to construct harmonic measures in the complex plane, see \cite{Kenig-Zhao} by Kenig and Zhao.  Another problem related to local smoothness of the Hilbert transform is the
local smoothness of the conjugate functions in the periodical case, and we refer the reader to the works by Dyachenko in \cite{Dyach_2} and \cite{Dyach_3} for recent progress on that front. 
We also remark that in this paper we assume $f \in L^p(\mathbb{R}),$ where $ 1\leq p < \infty$, however the same results would hold in the case $p= \infty$  when a modified Hilbert transform is used.

\vspace{.2in}

\section{The Operator $T$ and its H\"older Continuity}
\label{section:study_T}

\vspace{.1in}


\begin{defn}
	\label{eqn:transform_defn}
	Given $f \in C^{0, \alpha}([a,b])$, define the function $Tf$ on $[a,b]$ by 
	\begin{equation}
		\label{eqn:defn_T(f)}
		Tf (x) \,\, = \,\, \, \int _a^b  \frac{f(y)-f(x)}{y-x} \, dy.
	\end{equation}
\end{defn}

\noindent




\vspace{.1in}
\noindent
It is straightforward to verify that $Tf$ is well-defined for all H\"older continuous functions, but if $f$ is only continuous, then $Tf$ may not be well defined.  For example, let



 \begin{equation*}
		f(x) = 
		\begin{cases}
			\frac{1}{\ln x} & \,\, \text{if} \,\,\,\,  0< x \leq \frac{1}{2},  \\
			\noalign{\medskip}
			0  & \,\, \text{if} \,\,\,\,  x \leq 0.
		\end{cases}
	\end{equation*}
Then $f$ is continuous but not H\"older continuous, and $Tf(0)$ is not well defined.  
In addition, let $\eta(x)$ be a function that is compactly supported on $[-\frac{1}{2}, \frac{1}{2}]$ and $\eta =1$ on $[-\frac{1}{4}, \frac{1}{4}]$.  By definition $f\eta$ is $C^0$ but not $C^1$, and it is straightforward to check that its Hilbert transform $H(f\eta)$ is undefined at $0$.  Thus if a function is $C^0$ but not $C^1$, then its Hilbert transform may not be $C^0$.

\vspace{.1in}

In this section we will study various properties of the operator $T$ that eventually lead to the H\"older continuity of the derivatives of $T$, a key ingredient in the proof of Theorem \ref{thm:main}.   The first result is an estimate of the H\"older semi-norm of $T$.

\vspace{.1in}

\begin{thm}
	\label{thm:holder}
	If $f \in C^{0,\alpha}\left ( [a,b] \right )$ for some  $0<\alpha <1$, then for any $x_1, x_2 \in [a,b]$ ,
	\begin{equation}
		\label{eqn:Tf_difference}
		\left | Tf(x_2) - Tf(x_1) \right |  
		\leq   C \| f \|_{C^{0, \alpha}} |x_2 - x_1|^{\alpha} \left ( \ln \left \vert \frac{b-a}{x_2 - x_1 } \right \vert  +1 \right ).
	\end{equation}
In particular, $Tf \in C^{0, \alpha_0}\left ( [a,b] \right ) $ for any $0< \alpha_0 < \alpha$, and
	the map $T:  C^{0,\alpha} \to C^{0, \alpha_0}$ is a bounded operator.
 	
\end{thm}

\noindent
\textbf{Remark}: Theorem \ref{thm:holder} is a special case of Theorem \ref{thm:T_k_Holder} below. 
Both estimates in Theorem \ref{thm:holder} and Theorem \ref{thm:T_k_Holder} are sharp, i.e. the results cannot be improved to $ \alpha_0=\alpha$.  A simple counter-example is $f(x)=x^{\alpha}$ and detailed calculations are given in the proof of Lemma \ref{lemma:x^alpha_finite_Hilbert} in Section \ref{section:example}; we will discuss this in the remark immediately following the proof of Lemma \ref{lemma:x^alpha_finite_Hilbert}.

\vspace{.1in}
\noindent
\textit{Proof of Theorem \ref{thm:holder}}: 
Without loss of generality we can assume $a=0$, $b=1$, and $x_1 < x_2$.
Denote $$Q_{1}(y)=\frac{f(y)-f(x_1)}{y-x_1}  \hspace{.2in} \text{and} \hspace{.2in} Q_{2}(y)=\frac{f(y)-f(x_2)}{y-x_2}. $$ 
Since $f$ is $C^{0, \alpha}$, we have
\begin{equation}
	\label{eqn:Q_1_estimate}
	|Q_1(y)| \leq \| f \|_{C^{0, \alpha}} |y-x_1|^{\alpha -1}
\end{equation}
and
\begin{equation}
	\label{eqn:Q_2_estimate}
	|Q_2(y)| \leq \| f \|_{C^{0, \alpha}} |y-x_2|^{\alpha -1}.
\end{equation}

\noindent
We write
\begin{eqnarray}
	\label{eqn:difference_3_terms}
& &	\left | Tf(x_1) - Tf(x_2) \right | \nonumber   \\
& = & \left |  \int_0^1 \left ( Q_1(y) -Q_2(y) \right ) dy \right | \nonumber \\
\hspace{.3in} 	& \leq  &  \left |  \int_0^{x_1} \left ( Q_1(y) -Q_2(y) \right )dy \right | + \left | \int_{x_1}^{x_2} \left ( Q_1(y) -Q_2(y) \right )dy \right | + \left | \int_{x_2}^1 \left ( Q_1(y) -Q_2(y) \right )dy \right | \nonumber \\
	& := & I_1 + I_2 + I_3. 
\end{eqnarray}

 \vspace{.1in}
 \noindent
 We first estimate $$\displaystyle I_2=\left | \int_{x_1}^{x_2} \left ( Q_1(y) -Q_2(y) \right )dy \right | .
 $$
 \noindent
 By (\ref{eqn:Q_1_estimate}) and (\ref{eqn:Q_2_estimate}),
 
 \begin{eqnarray}
 	\label{eqn:int_x1_to_x2}
 	I_2 & \leq &  \int_{x_1}^{x_2} \left | Q_1(y) \right | dy +  \int_{x_1}^{x_2} \left | Q_2(y) \right | dy \nonumber \\
 	& \leq & \| f \|_{C^{0, \alpha}}\int_{x_1}^{x_2} |y-x_1|^{\alpha -1} dy + \| f \|_{C^{0, \alpha}} \int_{x_1}^{x_2} |y-x_2|^{\alpha -1} dy \nonumber \\
 	& = &  \frac{2\| f \|_{C^{0, \alpha}}}{\alpha}  (x_2-x_1)^{\alpha}.
 \end{eqnarray}

\vspace{.1in}
\noindent 
Next, we estimate $$\displaystyle I_1=\left | \int_0^{x_1} \left ( Q_1(y) -Q_2(y) \right )dy \right | .
$$
 	
 	\vspace{.1in}
 	 	\noindent
 	\textbf{Case 1}: $x_1 < x_2 - x_1$, which implies $x_2 < 2 (x_2-x_1)$.
 	 \vspace{.1in}
 	\noindent
 	Then by (\ref{eqn:Q_1_estimate}) and (\ref{eqn:Q_2_estimate})
 	\begin{equation*}
 		 \int_0^{x_1} \left | Q_{1}(y) \right | dy  \leq \int_0^{x_1} \| f \|_{C^{0, \alpha}} |y-x_1|^{\alpha -1} dy = \frac{\| f \|_{C^{0, \alpha}}}{\alpha} x_1^{\alpha} < \frac{ \| f \|_{C^{0, \alpha}}}{\alpha} (x_2 - x_1)^{\alpha}
 	\end{equation*}
 and
 \allowdisplaybreaks
 	\begin{eqnarray*}
 		 \int_0^{x_1} \left | Q_2 (y) \right | dy & \leq & \int_0^{x_1}  \| f \|_{C^{0, \alpha}}  |y-x_2|^{\alpha -1} dy \\
 		 & = & \frac{ \| f \|_{C^{0, \alpha}}}{\alpha} \left (  x_2^{\alpha} - (x_2 - x_1)^{\alpha}  \right ) \\
 		 & < &  \frac{ \| f \|_{C^{0, \alpha}}}{\alpha}x_2^{\alpha}  \\
 		 & < & \frac{2^{\alpha}\| f \|_{C^{0, \alpha}} }{\alpha}(x_2 - x_1)^{\alpha} .
 	\end{eqnarray*}
 
 \vspace{.05in}
 \noindent
 Therefore, $\displaystyle \left |  I_1 \right | \leq C \| f \|_{C^{0, \alpha}}(x_2 - x_1)^{\alpha} $, where $C$ only depends on $\alpha$.
 
 	 \vspace{.1in}
 	 	\noindent
 	 \textbf{Case 2}: $x_1 \geq x_2 - x_1$.  Denote	 
 \begin{equation}
 	\label{eqn:defn_h}
 	h= \frac{x_2 - x_1}{2}.
 \end{equation}
 We will integrate $\displaystyle \int_0^{x_1-h} \left ( Q_1(y) -Q_2(y) \right )dy$ \, and \, $\displaystyle \int_{x_1-h}^{x_1} \left ( Q_1(y) -Q_2(y) \right )dy$ \, separately.

\noindent
Note that 
\begin{eqnarray*}
	Q_1 - Q_2 & = & \frac{f(y)-f(x_1)}{y-x_1} - \frac{f(y)-f(x_2)}{y-x_2} \\
	& = & \frac{f(x_2)-f(x_1)}{y-x_1} + \left [ \frac{1}{y-x_1} - \frac{1}{y-x_2}   \right ] \left ( f(y)-f(x_2) \right ).
\end{eqnarray*}

\noindent
Since $f$ is $C^{0, \alpha}$, we have 

\begin{eqnarray}
	\label{eqn:int_to_x1-h_1st_term}
	\int_0^{x_1-h} \left | \frac{f(x_2)-f(x_1)}{y-x_1} \right | dy & \leq & \| f \|_{C^{0, \alpha}} |x_2-x_1|^{\alpha} \int_0^{x_1-h}   \frac{1}{x_1 - y}  dy \nonumber \\
	& = & \| f \|_{C^{0, \alpha}}|x_2-x_1|^{\alpha} \left (-\ln h + \ln x_1 \right ) \nonumber \\
	& < &  \| f \|_{C^{0, \alpha}} |x_2-x_1|^{\alpha}  \ln \left \vert \frac{2}{x_2 - x_1} \right \vert,
\end{eqnarray}
where in the last step we have used $\ln x_1 <0$. To estimate the integral of $$\left [ \frac{1}{y-x_1} - \frac{1}{y-x_2}   \right ] \left ( f(y)-f(x_2) \right ) $$ on $0 \leq y \leq x_1-h$, let $$v=\frac{x_1 - y}{x_2 - x_1}, $$ then  $$ \frac{1}{2} \leq v \leq \frac{x_1}{x_2-x_1} , \hspace{.2in}  x_1 - y = v(x_2 - x_1), \hspace{.2in} \text{and} \hspace{.2in} x_2 - y = (v+1)(x_2 - x_1). $$  Thus
\allowdisplaybreaks
\begin{eqnarray*}
	\left |  \left [ \frac{1}{y-x_1} - \frac{1}{y-x_2}   \right ] \left ( f(y)-f(x_2) \right )     \right | & = & \frac{1}{v(v+1)} \left | \frac{f(y)-f(x_2)}{x_2 - x_1}   \right | \\
	& \leq & \frac{1}{v(v+1)} \cdot \frac{\| f \|_{C^{0, \alpha}} |y - x_2|^{\alpha} }{|x_2 - x_1|}  \\
	& = & \frac{1}{v(v+1)} \cdot \frac{\| f \|_{C^{0, \alpha}}(v+1)^{\alpha}| x_2 - x_1|^{\alpha} }{|x_2 - x_1|} \\
	& = & \frac{\| f \|_{C^{0, \alpha}} }{v(v+1)^{1-\alpha}|x_2 - x_1|^{1-\alpha}}.
\end{eqnarray*}

\noindent
Consequently,
\allowdisplaybreaks
\begin{eqnarray}
		\label{eqn:int_to_x1-h_2nd_term}
	& &	\int_0^{x_1-h} 	\left |  \left [ \frac{1}{y-x_1} - \frac{1}{y-x_2}   \right ] \left ( f(y)-f(x_2) \right )     \right | dy \nonumber \\
	& \leq & \int ^{\frac{x_1}{x_2 - x_1}}_{\frac{1}{2}} \frac{\| f \|_{C^{0, \alpha}} }{v(v+1)^{1-\alpha}|x_2 - x_1|^{1-\alpha}} (x_2 - x_1) \, dv \nonumber \\
		& = & \| f \|_{C^{0, \alpha}}|x_2 - x_1|^{\alpha} \int_{\frac{1}{2}} ^ {\frac{x_1}{x_2 - x_1}} \frac{1}{v(v+1)^{1-\alpha} } \, dv \nonumber \\
		& < &  \| f \|_{C^{0, \alpha}}|x_2 - x_1|^{\alpha} \int_{\frac{1}{2}} ^ {\infty} \frac{1}{v(v+1)^{1-\alpha} } \, dv \nonumber \\
		& = & C\| f \|_{C^{0, \alpha}} |x_2 - x_1|^{\alpha},
\end{eqnarray}
where $C =\displaystyle \int_{\frac{1}{2}} ^ {\infty} \frac{1}{v(v+1)^{1-\alpha} } \, dv $ only depends on $\alpha$.  Combining (\ref{eqn:int_to_x1-h_1st_term}) and (\ref{eqn:int_to_x1-h_2nd_term}), we know
\begin{eqnarray}
	\label{eqn:int_0_to_x1-h}
\hspace{.4in}	\int _0^{x_1-h} |Q_1(y) - Q_2(y)| dy & \leq &  C\| f \|_{C^{0, \alpha}}|x_2 - x_1|^{\alpha} +  \| f \|_{C^{0, \alpha}}|x_2-x_1|^{\alpha}  \ln \left \vert \frac{2}{x_2 - x_1 } \right \vert  .
\end{eqnarray}

\vspace{.1in}
\noindent
By (\ref{eqn:Q_1_estimate}) and (\ref{eqn:Q_2_estimate}),
\allowdisplaybreaks
\begin{eqnarray}
	\label{eqn:int_x1-h_to_x1}
	\int _{x_1 -h}^{x_1} \left | Q_1(y) - Q_2(y)  \right | dy & \leq & \| f \|_{C^{0, \alpha}} \int _{x_1 -h}^{x_1} |y-x_1|^{\alpha -1} dy + \| f \|_{C^{0, \alpha}} \int _{x_1 -h}^{x_1} |y-x_2|^{\alpha -1} dy \nonumber \\
	& = & \frac{\| f \|_{C^{0, \alpha}}}{\alpha} \left \vert  \frac{x_2 - x_1}{2} \right \vert ^{\alpha} + \| f \|_{C^{0, \alpha}}\frac{\left (  \frac{3}{2} \right )^{\alpha} -1}{\alpha} \left \vert x_2 - x_1 \right \vert ^{\alpha} \nonumber \\
	& = & C\| f \|_{C^{0, \alpha}} \left | x_2 - x_1 \right |^{\alpha} ,
\end{eqnarray}

\vspace{.1in}
\noindent
where $C$ only depends on $\alpha$. It then follows from (\ref{eqn:int_0_to_x1-h}) and (\ref{eqn:int_x1-h_to_x1}) that
\begin{equation}
	\label{eqn:int_0_to_x1}
	I_1 \leq C \| f \|_{C^{0, \alpha}} \left | x_2 - x_1 \right |^{\alpha} \left ( 1+   \ln \left \vert \frac{1}{x_2 - x_1 } \right \vert \right ),
\end{equation}
where $C$ only depends on $\alpha$.

\vspace{.1in}
\noindent
Lastly, we estimate $$\displaystyle I_3=\left | \int_{x_2}^1 \left ( Q_1(y) -Q_2(y) \right )dy \right | .
$$ As in the estimate for $I_1$ we can write 
$$	Q_1 - Q_2   =  \frac{f(x_2)-f(x_1)}{y-x_1} + \left [ \frac{1}{y-x_1} - \frac{1}{y-x_2}   \right ] \left ( f(y)-f(x_2) \right ),$$ and 
\begin{eqnarray}
	\label{eqn:int_x2_to_1_1st_term}
	\int_{x_2}^1 \left | \frac{f(x_2)-f(x_1)}{y-x_1} \right | dy & \leq & \| f \|_{C^{0, \alpha}} |x_2-x_1|^{\alpha} \int_{x_2}^1   \frac{1}{y - x_1 }  dy \nonumber \\
	& = & \| f \|_{C^{0, \alpha}} |x_2-x_1|^{\alpha} \big (\ln \vert 1-x_1 \vert - \ln \vert x_2 - x_1 \vert \big ) \nonumber \\
	& < &  \| f \|_{C^{0, \alpha}} |x_2-x_1|^{\alpha}  \ln \left \vert \frac{1}{x_2 - x_1} \right \vert,
\end{eqnarray}
where we have used $\ln(1-x_1)<0$.  To estimate the second term, let $$v=\frac{ y - x_1}{x_2 - x_1}, $$ then $1<v<\frac{ 1 - x_1}{x_2 - x_1} $ when $x_2 <y <1$, in addition, $$ y - x_1  = v(x_2 - x_1) \hspace{.2in} \text{and} \hspace{.2in} y - x_2  = (v-1)(x_2 - x_1). $$  Thus
\allowdisplaybreaks
\begin{eqnarray*}
	\left |  \left [ \frac{1}{y-x_1} - \frac{1}{y-x_2}   \right ] \left ( f(y)-f(x_2) \right )     \right | & = & \frac{1}{v(v-1)} \left | \frac{f(y)-f(x_2)}{x_2 - x_1}   \right | \\
	& \leq & \frac{1}{v(v-1)} \cdot \frac{\| f \|_{C^{0, \alpha}}|y - x_2|^{\alpha} }{|x_2 - x_1|}  \\
	& = & \frac{1}{v(v-1)} \cdot \frac{\| f \|_{C^{0, \alpha}}(v-1)^{\alpha}| x_2 - x_1|^{\alpha} }{|x_2 - x_1|} \\
	& = & \frac{\| f \|_{C^{0, \alpha}} }{v(v-1)^{1-\alpha}|x_2 - x_1|^{1-\alpha}}.
\end{eqnarray*}

\noindent
Consequently,
\allowdisplaybreaks
\begin{eqnarray}
	\label{eqn:int_x2_to_1_2nd_term}
	& &	\int_{x_2}^1 	\left |  \left [ \frac{1}{y-x_1} - \frac{1}{y-x_2}   \right ] \left ( f(y)-f(x_2) \right )     \right | dy \nonumber \\
	& \leq & \int _1^{\frac{1-x_1}{x_2 - x_1}} \frac{\| f \|_{C^{0, \alpha}} }{v(v-1)^{1-\alpha}|x_2 - x_1|^{1-\alpha}} (x_2 - x_1) \, dv \nonumber \\
	& = & \| f \|_{C^{0, \alpha}}|x_2 - x_1|^{\alpha} \int_1^{\frac{1- x_1}{x_2 - x_1}} \frac{1}{v(v-1)^{1-\alpha} } \, dv \nonumber \\
	& < &  \| f \|_{C^{0, \alpha}}|x_2 - x_1|^{\alpha} \int_1^ {\infty} \frac{1}{v(v-1)^{1-\alpha} } \, dv \nonumber \\
	& = & C\| f \|_{C^{0, \alpha}}|x_2 - x_1|^{\alpha},
\end{eqnarray}
where $C=\displaystyle \int_1^{\infty} \frac{1}{v(v-1)^{1-\alpha} } \, dv $ only depends on $\alpha$.  Combining (\ref{eqn:int_x2_to_1_1st_term}) and (\ref{eqn:int_x2_to_1_2nd_term}), we know
\begin{eqnarray}
	\label{eqn:int_x2_to_1}
	I_3 & \leq & \int _{x_2}^1 |Q_1(y) - Q_2(y)| dy \nonumber \\
	& \leq & C \| f \|_{C^{0, \alpha}} |x_2 - x_1|^{\alpha} \left ( 1+ \ln \left \vert \frac{1}{x_2 - x_1 } \right \vert  \right )  , 
\end{eqnarray}
where $C$ only depends on $\alpha$.

\vspace{.1in}
\noindent
Therefore, by (\ref{eqn:difference_3_terms}),  (\ref{eqn:int_x1_to_x2}), (\ref{eqn:int_0_to_x1}), and (\ref{eqn:int_x2_to_1}),
\begin{eqnarray*}
	 \left | Tf(x_1) - Tf(x_2) \right |  
	 & \leq &  C \| f \|_{C^{0, \alpha}} |x_2 - x_1|^{\alpha} \left ( 1+ \ln \left \vert \frac{1}{x_2 - x_1 } \right \vert  \right ) .
\end{eqnarray*}
For any $0 < \alpha_0 < \alpha$, 
\begin{eqnarray*}
  \frac{ \left | Tf(x_1) - Tf(x_2) \right | }{|x_1 - x_2|^{\alpha_0}} 
 & \leq & C \| f \|_{C^{0, \alpha }} |x_2 - x_1|^{\alpha  -\alpha _0 } \left ( 1 + \ln \left \vert \frac{1}{x_2 - x_1 } \right \vert \right )  .  
 \end{eqnarray*}
Hence $\| Tf \|_{ C^{0, \alpha_0}} \leq C\| f \|_{C^{0, \alpha}}$ and $ T:  C^{0,\alpha} \to C^{0, \alpha_0} $ is bounded.

\stop

\vspace{.2in}

In the rest of this section we will study the derivatives of $Tf$.  One important observation is an intriguing connection between the derivatives of the difference quotient and the remainder of the Taylor expansion of $f$.

\vspace{.1in}

\noindent
Recall that for any $f \in C^k([a,b])$,	
\begin{equation*}
	P_k(x,y) = \sum_{i =0}^{k} \dfrac{f^{(i)}(x)}{i!}(y-x)^i
\end{equation*}
is the $k$-th order Taylor polynomial of $f$ at $x$.

\vspace{.2in}

\begin{lemma}
	\label{lem:derivative_diff_quotient}
	For any $k \in \mathbb{N}$, if $f$ is a $ C^k$ function,  then $$ \frac{\partial ^k}{\partial x^k}\left ( \frac{f(y)-f(x)}{y-x} \right )  = k! \dfrac{f(y)-P_k(x,y)}{(y-x)^{k+1}}, $$ for $x \neq y$, where $P_k(x,y)$ is defined above.
\end{lemma}

\vspace{.1in}
\pf We will prove this lemma by induction.
When $k=1$, by direct calculations
$$ 	\frac{\partial }{\partial x} \left( \frac{f(y)-f(x)}{y-x} \right )  =  \dfrac{ f(y)-f(x)-f'(x)(y-x)}{(y-x)^2} = \dfrac{f(y)-P_1(x,y)}{(y-x)^{2}}. $$

\noindent
Suppose for $k \in \mathbb{N}$, all $C^k$ functions $f$ satisfy $$ \frac{\partial ^k}{\partial x^k} \left( \frac{f(y)-f(x)}{y-x} \right )  = k! \dfrac{f(y)-P_k(x,y)}{(y-x)^{k+1}}, $$ 

\vspace{.05in}
\noindent
then for any $f \in C^{k+1}$,

\allowdisplaybreaks
\begin{equation*}
	\frac{\partial ^{k+1} }{\partial x ^{k+1}} \left( \frac{f(y)-f(x)}{y-x} \right )   = k!  \frac{\partial }{\partial x} \left ( \dfrac{f(y)-P_k(x,y)}{(y-x)^{k+1}} \right ). 
\end{equation*}

\noindent
Since $$P_k(x,y) = f(x) + f'(x)(y-x) + \frac{f''(x)}{2}(y-x)^2 + \cdots + \frac{f^{(k)}(x)}{k!}(y-x)^k,  $$ we have
\begin{eqnarray*}
	\frac{\partial }{\partial x}   P_k(x,y) & = & f'(x) + \Big ( f''(x)(y-x) - f'(x) \Big )  + \left ( \frac{f'''(x)}{2}(y-x)^2 - f''(x) (y-x) \right ) \\
	&  & + \cdots + \left ( \frac{f^{(k+1)}(x)}{k!}(y-x)^k  - \frac{f^{(k)}(x)}{(k-1)!}(y-x)^{k-1} \right ) \\
	& = &  \frac{f^{(k+1)}(x)}{k!}(y-x)^k.
\end{eqnarray*}

\noindent
Thus 
\allowdisplaybreaks
\begin{eqnarray}
	\label{eqn:deriv_Taylor_remainder}
	\frac{\partial }{\partial x} \left ( \dfrac{f(y)-P_k(x,y)}{(y-x)^{k+1}} \right ) & = &  \dfrac{ - \left (\frac{\partial }{\partial x}   P_k(x,y) \right ) (y-x)^{k+1} +(k+1) \left ( f(y)-P_k(x,y) \right ) (y-x)^k }{(y-x)^{2k+2}} \nonumber \\
	& = & \dfrac{ - \frac{f^{(k+1)}(x)}{k!} (y-x)^{2k+1} +(k+1) \left ( f(y)-P_k(x,y) \right ) (y-x)^k }{(y-x)^{2k+2}} \nonumber \\
	& = & (k+1) \dfrac{ - \frac{f^{(k+1)}(x)}{(k+1)!} (y-x)^{k+1} +  f(y)-P_k(x,y)  }{(y-x)^{k+2}} \nonumber \\
	& = & (k+1) \dfrac{ f(y)-P_{k+1}(x,y)   }{(y-x)^{k+2}},
\end{eqnarray}
and consequently,
\allowdisplaybreaks
\begin{eqnarray*}
	\frac{\partial ^{k+1} }{\partial x ^{k+1}} \left( \frac{f(y)-f(x)}{y-x} \right )  
	& = & (k+1)! \dfrac{ f(y)-P_{k+1}(x,y)   }{(y-x)^{k+2}}.
\end{eqnarray*}
This completes the proof.
\stop

\vspace{.2in}

Lemma \ref{lem:derivative_diff_quotient} inspires us to define another operator $T_k$ as below, which will be used to study the derivatives of $Tf$. Its definition can be traced back to \cite{Liu-Pan-Zhang} and \cite{Pan_Zhang} in the complex setting.

\begin{defn}
	\label{defn:T_kf}
	For any $f \in C^{k, \alpha} \left ( [a,b] \right ) (0< \alpha <1)$, define
	\begin{equation}
		\label{eqn:T_k(f)}
		T_kf(x) = k! \int _a^b \frac{f(y)-P_k(x,y)}{(y-x)^{k+1}} dy,
	\end{equation}
	where $P_k(x,y)$ is the $k$-th degree Taylor polynomial of $f$ at $x$. 
\end{defn}

\vspace{.1in}
\noindent
Note that
\begin{eqnarray*}
	f(y)-P_{k}(x,y) & = &  f(y) -P_{k-1}(x,y) - \frac{f^{(k)}(x)}{k!}(y-x)^{k} \\
	& = &  \frac{f^{(k)}(\zeta)}{k!}(y-x)^{k} - \frac{f^{(k)}(x)}{k!}(y-x)^{k}
\end{eqnarray*}
for some $\zeta$ between $x$ and $y$.  Since $f$ is $C^{k, \alpha}$, it follows that
$$ \left \vert   f^{(k)}(\zeta) - f^{(k)}(x)  \right \vert \leq C |\zeta - x|^{\alpha} \leq C |y - x|^{\alpha} . $$
Hence 
\begin{equation}
	\label{eqn:Taylor_remainder_quotient_estimate}
	 \left |  \frac{f(y)-P_{k}(x,y)}{(y-x)^{k+1}}  \right | \leq C|y-x|^{\alpha -1},
\end{equation}

\vspace{.05in}
\noindent
which implies that $T_kf$ is well-defined.  Later in this section we will prove the H\"older continuity of $T_k$ in Theorem \ref{thm:T_k_Holder}, but before doing so we need a few preparatory lemmas.

\vspace{.1in}

\begin{lemma}
	\label{lem:Talor_expansion_Holer}
	Suppose a function $f \in C^{k+1, \alpha} \left ( [a,b] \right ) (0< \alpha <1)$, then the function $$ \frac{f(y)-P_k(x,y)}{(y-x)^{k+1}} $$ is well-defined on $[a,b] \times [a,b]$ and H\"older continuous with exponent $\alpha$ in $x$ and $y$.
\end{lemma}

\vspace{.1in}
\pf 
Since $f^{(k+1)}$ is continuous, the Taylor expansion of $f$ at $x$ takes the form 
\begin{eqnarray*}
	f(y) & = &  P_k(x,y) + \frac{1}{k!} \int _x ^y (y-t)^k f^{(k+1)} (t) \, dt \\
	& = & P_k(x,y) +  \frac{(y-x)^{k+1}}{k!} \int _0 ^1 (1-s)^k f^{(k+1)} \big (sy+(1-s)x\big ) \, ds  ,
\end{eqnarray*}

\vspace{.05in}
\noindent
after a change of variable $s=\frac{t-x}{y-x}$.  Thus $$ \dfrac{f(y)-P_k(x,y)}{(y-x)^{k+1}} = \frac{1}{k!} \int _0 ^1 (1-s)^k f^{(k+1)} \big (sy+(1-s)x\big ) \, ds $$ is well defined in $[a,b] \times [a,b]$.
  
\noindent
For any $a \leq x_1,x_2,y_1,y_2 \leq b$ and $(x_1, y_1) \neq (x_2, y_2)$, 
\begin{eqnarray*}
& &	 \left | \int _0 ^1  (1-s)^k f^{(k+1)}\big (sy_1+(1-s)x_1 \big ) \, ds -  \int _0 ^1 (1-s)^k f^{(k+1)}\big (sy_2 +(1-s)x_2 \big ) \, ds \right | \\
	& \leq & \int _0^1 \left |  f^{(k+1)}\big (sy_1+(1-s)x_1 \big ) - f^{(k+1)}\big (sy_2+(1-s)x_2\big )  \right | ds \\
	& \leq &  \|f\|_{C^{k+1, \alpha}} \left ( \left |  x_1 - x_2  \right |^{\alpha} + \left |  y_1 - y_2  \right |^{\alpha} \right ).
\end{eqnarray*}
Thus $\displaystyle \frac{f(y)-P_k(x,y)}{(y-x)^{k+1}}$ is H\"older continuous with exponent $\alpha$ in $x$ and $y$. 

\stop

\vspace{.2in}
Next, we will prove that the differentiation and integral signs are interchangeable for the operator $T$.  This turns out to be a subtle issue, and the following folklore elementary lemma will play a crucial role in what follows.

\vspace{.1in}

\begin{lemma}
	\label{lem:deriv_convergence}
	Let $\{ f_n (x) \}$ be a sequence of differentiable functions from $[a,b] $ to $\mathbb{R}$. Suppose
	$f_n \to f$ uniformly on $[a,b]$ and
	$f'_n \to g$ uniformly on $[a,b]$, where $f,g: [a,b] \to \mathbb{R}$.
	Then $f'(x)=g(x)$ for every $x \in [a,b]$.
\end{lemma}

\noindent
\textbf{Remark:} The key point of this lemma is that $f'=g$ at every point in $[a,b]$.  It is easy to show $f'=g$ almost everywhere, since  $$\int_a ^x f'(t) dt = f(x) - f(a) = \lim_{n \to \infty} f_n(x) - f(a) = \lim_{n \to \infty} \int_a ^x f_n'(t) dt= \int_a ^x g(t) dt.$$  If $g$ is continuous, this would imply $f'=g$.  However, $g$ could be discontinuous in general.  For example, define a trivial sequence
\begin{equation*}
	f_n(x) = f(x) =
	\begin{cases}
		x^2 \sin \left ( \frac{1}{x^2} \right ) & \,\, \text{if} \,\,\,\,   x \neq 0,  \\
		\noalign{\medskip}
		0  & \,\, \text{if} \,\,\,\,x=0  
	\end{cases}
\end{equation*} 
for all $n\in \mathbb{N}$, then $f_n$ is differentiable, but
\begin{equation*}
	f'_n(x) = g(x) =
	\begin{cases}
		2x \sin \left ( \frac{1}{x^2} \right ) - \frac{2}{x} \cos \left ( \frac{1}{x^2} \right ) & \,\, \text{if} \,\,\,\,   x \neq 0,  \\
		\noalign{\medskip}
		0  & \,\, \text{if} \,\,\,\,x=0  
	\end{cases}
\end{equation*} 
is not continuous.
  The proof of this lemma is not trivial, so it is included in Appendix A.

\stop

\vspace{.1in}

\begin{lemma}
	\label{lem:deriv_integral_Taylor_remainder}
	 Suppose a function $f \in C^{k+1, \alpha} \left ( [a,b] \right )$, then
	\begin{equation}
		\label{eqn:deriv_integral_Taylor_remainder}
		\frac{d}{dx} \left ( \int _a^b \frac{f(y)-P_k(x,y)}{(y-x)^{k+1}} dy \right ) = \int _a^b \frac{\partial }{\partial x} \left ( \frac{f(y)-P_k(x,y)}{(y-x)^{k+1}} \right ) dy. \nonumber
	\end{equation} 
\end{lemma}

\pf
Denote $$F(x) = \int _a^b \frac{f(y)-P_k(x,y)}{(y-x)^{k+1}} \, dy.$$  In order to avoid the singularity in the integral, we take any sequence of positive numbers $  \epsilon_n  \to 0 $ and define $$F_n (x) = \int_a ^{x-\epsilon_n} \frac{f(y)-P_k(x,y)}{(y-x)^{k+1}} \, dy + \int_{x+\epsilon_n} ^b \frac{f(y)-P_k(x,y)}{(y-x)^{k+1}} \, dy. $$  Then since $f  \in C^{k+1, \alpha} \left ( [a,b] \right ), $ it follows by the Taylor expansion $\displaystyle \left \vert \frac{f(y)-P_k(x,y)}{(y-x)^{k+1}} \right \vert \leq C$.  Consequently, 
\begin{equation*}
	\left | F(x) - F_n(x)   \right | \,\, = \,\, \left \vert  \int_{x-\epsilon_n}^{x+\epsilon_n} \frac{f(y)-P_k(x,y)}{(y-x)^{k+1}} \, dy \right \vert \,\, \leq  \,\, C\epsilon_n.
\end{equation*}
Thus $F_n \to F$ uniformly.  

\noindent
Next, we consider the limit of the derivative of $F_n$.  By the Leibniz formula

\allowdisplaybreaks
\begin{eqnarray*}
	F'_n(x) & = & \int_a ^{x-\epsilon_n} \frac{\partial }{\partial x} \left ( \frac{f(y)-P_k(x,y)}{(y-x)^{k+1}} \right ) dy + \frac{f(x-\epsilon_n)-P_k(x,x-\epsilon_n)}{(-\epsilon_n)^{k+1}} \\
	& + & \int_{x+\epsilon_n} ^b \frac{\partial }{\partial x} \left ( \frac{f(y)-P_k(x,y)}{(y-x)^{k+1}} \right ) dy - \frac{f(x+\epsilon_n)-P_k(x,x+\epsilon_n)}{\epsilon_n^{k+1}}.
\end{eqnarray*}
It follows that 
\begin{eqnarray*}
& & \left | F'_n(x) -  \int _a^b \frac{\partial }{\partial x} \left ( \frac{f(y)-P_k(x,y)}{(y-x)^{k+1}} \right ) dy \right | \\
& \leq  & \left | \int_{x-\epsilon_n}^{x+\epsilon_n} \frac{\partial }{\partial x} \left ( \frac{f(y)-P_k(x,y)}{(y-x)^{k+1}} \right ) dy \right |
 +   \left | \frac{f(x-\epsilon_n)-P_k(x,x-\epsilon_n)}{(-\epsilon_n)^{k+1}} - \frac{f(x+\epsilon_n)-P_k(x,x+\epsilon_n)}{\epsilon_n^{k+1}}\right |.
\end{eqnarray*}

\vspace{.1in}
\noindent
By (\ref{eqn:deriv_Taylor_remainder}),
\begin{eqnarray*}
	 \left | \int_{x-\epsilon_n}^{x+\epsilon_n} \frac{\partial }{\partial x} \left ( \frac{f(y)-P_k(x,y)}{(y-x)^{k+1}} \right ) dy \right | & = & \left | \int_{x-\epsilon_n}^{x+\epsilon_n}  (k+1) \dfrac{ f(y)-P_{k+1}(x,y)   }{(y-x)^{k+2}} dy \right |.
\end{eqnarray*}
By (\ref{eqn:Taylor_remainder_quotient_estimate}) with $k$ replaced by $k+1$, we know
$$ \left |  \frac{f(y)-P_{k+1}(x,y)}{(y-x)^{k+2}}  \right | \leq C|y-x|^{\alpha -1}.$$
Thus
$$  \left | \int_{x-\epsilon_n}^{x+\epsilon_n} \frac{\partial }{\partial x} \left ( \frac{f(y)-P_k(x,y)}{(y-x)^{k+1}} \right ) dy \right | \leq C \int_{x-\epsilon_n}^{x+\epsilon_n}  |y-x|^{\alpha -1} dy = C \epsilon_n ^ {\alpha}. $$

\vspace{.1in}
\noindent
By Lemma \ref{lem:Talor_expansion_Holer}, $\displaystyle \frac{f(y)-P_k(x,y)}{(x-y)^{k+1}}$ is H\"older in $y$, so
$$  \left | \frac{f(x-\epsilon_n)-P_k(x,x-\epsilon_n)}{(-\epsilon_n)^{k+1}} - \frac{f(x+\epsilon_n)-P_k(x,x+\epsilon_n)}{\epsilon_n^{k+1}}\right | \leq C \left | ( x- \epsilon_n ) - (x+ \epsilon_n) \right |^{\alpha} \leq C\epsilon_n ^{\alpha}. $$

\noindent
Therefore $$ \left | F'_n(x) -  \int _a^b \frac{\partial }{\partial x} \left ( \frac{f(y)-P_k(x,y)}{(y-x)^{k+1}} \right ) dy \right |  \leq C\epsilon_n ^{\alpha},$$
which implies that $$\displaystyle F'_n \to \int _a^b \frac{\partial }{\partial x} \left ( \frac{f(y)-P_k(x,y)}{(y-x)^{k+1}} \right ) dy$$ uniformly.  Then by Lemma \ref{lem:deriv_convergence} we know $$ F'(x) = \int _a^b \frac{\partial }{\partial x} \left ( \frac{f(y)-P_k(x,y)}{(y-x)^{k+1}} \right ) dy, $$ and this completes the proof of this lemma.

\stop

\vspace{.2in}

\noindent
\textbf{Remark}:
	It is worth noting that Lemma \ref{lem:deriv_integral_Taylor_remainder} will not be true if $f$ is only $C^{k+1}$ rather than $C^{k+1, \alpha}$, as shown in the next example.

\vspace{.1in}

\noindent
\textbf{Example}:
	 Let \begin{equation*}
		f(x) = 
		\begin{cases}
			\frac{x}{\ln x } & \,\, \text{if} \,\,\,\,  0< x \leq \frac{1}{2},  \\
			\noalign{\medskip}
			0  & \,\, \text{if} \,\,\,\,  x = 0.
		\end{cases}
	\end{equation*}
	$f \in C^1\left ([0, \frac{1}{2}] \right )$, but it is not $C^{1,\alpha}$ at $0$ for any $\alpha \in (0,1)$. 
	By definition
	\begin{eqnarray*}
		\frac{\partial }{\partial x} \left ( \frac{f(y)-f(x)}{y-x} \right ) \Bigg \vert _{x=0} 
		 & = & \lim_{h \to 0} \dfrac{\frac{f(y)-f(h)}{y-h}  - \frac{f(y)}{y}  }{h} \\
		 & = & \lim_{h \to 0} \dfrac{\frac{\frac{y}{\ln y } - \frac{h}{\ln h } }{y-h}  - \frac{\frac{y}{\ln y }}{y}  }{h} \\
		& = &  \lim_{h \to 0} \dfrac{\frac{1}{\ln y} - \frac{1}{\ln h} }{y-h} \\
		& = & \frac{1}{y\ln y}.
	\end{eqnarray*}

\noindent
Since $\displaystyle \int_0^{\frac{1}{2}} \frac{1}{y\ln y} \, dy = - \infty,$ it follows that $\displaystyle \int_0^{\frac{1}{2}} 	\frac{\partial }{\partial x} \left ( \frac{f(y)-f(x)}{y-x} \right )  \, dy$ is not defined at $x=0$.

\vspace{.2in}
The following is an application of Lemma \ref{lem:deriv_integral_Taylor_remainder} to higher derivatives of $Tf$.

\vspace{.1in}

\begin{thm}
	\label{thm:switch_deriv_integral}
	 Suppose a function $f \in C^{k, \alpha} \left ( [a,b] \right ) (0< \alpha <1)$, then
	\begin{equation}
		\label{eqn:deriv_integral_change}
		\frac{d^k}{dx^k} \int_a^b \frac{f(y)-f(x)}{y-x} dy = \int _a^b \frac{\partial ^k}{\partial x^k} \left ( \frac{f(y)-f(x)}{y-x} \right ) dy. \nonumber
	\end{equation} 
\end{thm}

\pf
If $f \in C^{1, \alpha} \left ( [a,b] \right )$, by Lemma \ref{lem:deriv_integral_Taylor_remainder}
$$ 	\frac{d}{dx} \int_a^b \frac{f(y)-f(x)}{y-x} dy = \int _a^b \frac{\partial }{\partial x} \left ( \frac{f(y)-f(x)}{y-x} \right ) dy,$$ so the theorem is true when $k=1$.  Now suppose $$	\frac{d^k}{dx^k} \int_a^b \frac{f(y)-f(x)}{y-x} dy = \int _a^b \frac{\partial ^k}{\partial x^k} \left ( \frac{f(y)-f(x)}{y-x} \right ) dy$$ for all $f \in C^{k, \alpha} \left ( [a,b] \right ) $,  we will show the theorem is true in $k+1$ for all $f \in C^{k+1, \alpha} \left ( [a,b] \right ) $.  Since

\allowdisplaybreaks
\begin{eqnarray*}
\frac{d^{k+1}}{dx^{k+1}} \int_a^b \frac{f(y)-f(x)}{y-x} dy & = & \frac{d}{dx} \left ( \frac{d^k}{dx^k} \int_a^b \frac{f(y)-f(x)}{y-x} dy   \right )	\\
& = &  \frac{d}{dx} \left ( \int _a^b \frac{\partial ^k}{\partial x^k} \left ( \frac{f(y)-f(x)}{y-x} \right ) dy  \right )	\\
& = & \frac{d}{dx} \left ( \int _a^b k! \dfrac{f(y)-P_k(x,y)}{(y-x)^{k+1}} dy  \right ) \hspace{.2in} \text{(by Lemma \ref{lem:derivative_diff_quotient})} \\
& = &  \int _a^b \frac{\partial }{\partial x} \left ( k! \frac{f(y)-P_k(x,y)}{(y-x)^{k+1}} \right ) dy \hspace{.2in} \text{(by Lemma \ref{lem:deriv_integral_Taylor_remainder})} \\
& = & \int_a^b  \frac{\partial ^{k+1}}{\partial x^{k+1}} \left ( \frac{f(y)-f(x)}{y-x} \right ) dy. \hspace{.2in} \text{(by Lemma \ref{lem:derivative_diff_quotient})} 
\end{eqnarray*}
Thus the theorem is proved by induction.

\stop

\vspace{.1in}

\noindent
\textbf{Remark}:	By Theorem \ref{thm:switch_deriv_integral}, Lemma \ref{lem:derivative_diff_quotient}, and Definition \ref{defn:T_kf} we know
	\begin{equation}
		\label{eqn:Tf_kth_derivative}
		\frac{d^k}{dx^k} Tf(x) \,\,  = \,\, k!\int _a^b \frac{f(y)-P_k(x,y)}{(y-x)^{k+1}} dy \,\, = \,\,  T_kf(x) .
	\end{equation}

\vspace{.2in}

\begin{thm}
	\label{thm:T_k_Holder}
	 Suppose a function $f \in C^{k, \alpha} \left ( [a,b] \right )$ for $ 0< \alpha <1$ and $k \geq 0$. Then for any $x_1, x_2 \in [a,b]$,
	 \begin{equation}
	 	\label{eqn:T_kf_difference}
	 	\left | T_k f(x_2) - T_k f (x_1) \right |  
	 	\leq   C \| f \|_{C^{k, \alpha}} |x_2 - x_1|^{\alpha} \left ( \ln \left \vert \frac{b-a}{x_2 - x_1 } \right \vert +1 \right )  .
	 \end{equation}
In particular, $T_k f \in C^{0, \alpha_0} \left ( [a,b] \right ) $ for all $0< \alpha_0 < \alpha $ and the map $T_k:  C^{k,\alpha} \to C^{0, \alpha_0}$ is a bounded operator; furthermore, $Tf \in C^{k, \alpha_0} \left ( [a,b] \right ).$
\end{thm}

\pf We will prove (\ref{eqn:T_kf_difference}) by induction.  Note that (\ref{eqn:T_kf_difference}) is true when $k=0$ by Theorem \ref{thm:holder}.  Now suppose it is true for $k$, we will show that it is also true for $k+1$.

\vspace{.1in}
\noindent
For any $f \in C^{k+1, \alpha}([a,b])$, by  (\ref{eqn:Tf_kth_derivative}) and Theorem \ref{thm:switch_deriv_integral}.
\allowdisplaybreaks
\begin{eqnarray*}
	T_{k+1} f(x) & = & \frac{d^{k+1}}{dx^{k+1}} Tf(x) \\
	& = & \frac{d}{dx} \left ( \int _a^b \frac{\partial ^k}{\partial x^k} \left ( \frac{f(y)-f(x)}{y-x} \right ) dy \right ) \hspace{.5in} \text{(by Theorem \ref{thm:switch_deriv_integral})} \\
	& = & \frac{d}{dx} \left ( \int _a^b k! \dfrac{f(y)-P_k(x,y)}{(y-x)^{k+1}} \, dy \right ) \hspace{.7in} \text{(by Lemma \ref{lem:derivative_diff_quotient})}\\
	& = & k!  \int _a^b \frac{\partial }{\partial x} \left ( \frac{f(y)-P_k(x,y)}{(y-x)^{k+1}} \right ) dy \hspace{.7in} \text{(by Lemma \ref{lem:deriv_integral_Taylor_remainder})}  \\
		& = & k!  \int _{a-x}^{b-x} \frac{\partial }{\partial x} \left ( \frac{f(x+z)-P_k(x,x+z)}{z^{k+1}} \right ) dz \hspace{.4in} (\text{let} \,\, z=y-x)  \\
	& 	= & k!  \int _{a-x}^{b-x}  \dfrac{f'(x+z)-\displaystyle \sum_{i=0}^k \frac{f^{(i+1)}(x)}{i!}z^i}{z^{k+1}} \, dz. \\
		& = & k!  \int_a^b \dfrac{f'(y)-\displaystyle \sum_{i=0}^k \frac{f^{(i+1)}(x)}{i!}(y-x)^i}{(y-x)^{k+1}}  dy \hspace{.7in} (\text{let} \,\, y=z+x)  .
\end{eqnarray*}

\vspace{.1in}
\noindent
Note that by (\ref{eqn:Tf_kth_derivative})
$$ k! \int_a^b \dfrac{f'(y)-\displaystyle \sum_{i=0}^k \frac{f^{(i+1)}(x)}{i!}(y-x)^i}{(y-x)^{k+1}}  dy = T_k f'(x), $$ hence
$$ T_{k+1}f =  T_k f'. $$
 Since $f' \in C^{k, \alpha} \left ( [a,b] \right )$, by the inductive hypothesis
 $$
 	\left | T_k f'(x_2) - T_k f' (x_1) \right |  
 	\leq   C \| f' \|_{C^{k, \alpha}} |x_2 - x_1|^{\alpha}  \left (\ln \left \vert \frac{b-a}{x_2 - x_1 } \right \vert +1 \right ) .
 $$ 
Therefore,
\begin{eqnarray*}
	\left | T_{k+1} f (x_2) - T_{k+1} f(x_1) \right |  
&	\leq &   C \| f' \|_{C^{k, \alpha}} |x_2 - x_1|^{\alpha}  \left (\ln \left \vert \frac{b-a}{x_2 - x_1 } \right \vert +1 \right ) \\
& \leq &   C \| f \|_{C^{k+1, \alpha}} |x_2 - x_1|^{\alpha}  \left (\ln \left \vert \frac{b-a}{x_2 - x_1 } \right \vert +1 \right ) .
\end{eqnarray*}

\noindent
Thus (\ref{eqn:T_kf_difference}) is true by induction.  It follows that $\| T_{k} f \| \in C^{0, \alpha_0} \left (  [a,b] \right ) $ for all $0< \alpha_0 < \alpha $ and  $T_{k}: C^{k, \alpha} \left ( [a,b] \right ) \to C^{0, \alpha_0} \left ( [a,b] \right ) $ is a bounded operator.  By (\ref {eqn:Tf_kth_derivative}) it also implies that $Tf \in C^{k, \alpha_0} \left ( [a,b] \right ).$
This completes the proof of this theorem.

\stop

\vspace{.2in}

\section{Proof of  Theorem \ref{thm:main}}
\label{section:main_proof}

\vspace{.1in}

Now we are in the position to prove Theorem \ref{thm:main}.
Since  $f \in L^p(\mathbb{R})$, by classical theory its Hilbert transform $Hf$ is defined a.e.
Suppose $f$ is $C^{k, \alpha} (k \geq 0)$ on an interval $[a, b]$, we will show $Hf$ is $C^{k, \alpha_0}  \left (0<\alpha_0 < \alpha \right )$ on any closed subinterval of $(a,b)$.  Letting $\chi_{(a,b)} $ be the characteristic function on $(a,b)$, we can write 
\begin{equation*}
	\label{eqn:Hf_decompose}
	Hf = H \left ( \chi_{(a,b)} f  \right ) + H \left (  (1- \chi_{(a,b)}) f \right )
\end{equation*}
and we will look at the two terms separately.

\vspace{.1in}
\noindent
For any $x \in [a_0, b_0] \subset (a,b)$, 
\begin{eqnarray}
	\label{eqn:H_chi_f}
	H \left ( \chi_{(a,b)} f \right )(x)
	& = &  p.v. \, \frac{1}{\pi}\int_a^b \frac{f(y)}{x-y} dy \nonumber \\
	& = & p.v. \, \frac{1}{\pi} \int_a^b \frac{f(y) - f(x)}{x-y} dy +   p.v. \, \frac{1}{\pi}  \int_a^b \frac{f(x)}{x-y} dy \nonumber \\
	& = & -\frac{1}{\pi}Tf(x)  + \frac{f(x)}{\pi}\ln \left \vert \frac{x-a}{x-b} \right \vert, 
\end{eqnarray}

\noindent
where we have used $$ p.v. \,  \int_a^b \frac{1}{x-y} dy = \ln \left \vert \frac{x-a}{x-b} \right \vert. $$ 

\vspace{.1in}
\noindent
For any $0<\alpha_0 < \alpha$, Theorem \ref{thm:T_k_Holder} implies $Tf \in C^{k, \alpha_0}\left ( [a_0, b_0] \right)$, so by (\ref{eqn:H_chi_f}) we know that
 $	H \left ( \chi_{(a,b)} f \right ) \in C^{k, \alpha_0}\left ( [a_0, b_0] \right)$.
By the following Lemma \ref{lemma:outside_analytic},  $H \left (  (1- \chi_{(a,b)}) f \right )$ is real analytic on  $ [a_0, b_0]$.  Therefore $Hf\in C^{k, \alpha_0}\left ( [a_0, b_0] \right)$, and this completes the proof of Theorem \ref{thm:main}.

\stop

\vspace{.1in}
\begin{lemma}
	\label{lemma:outside_analytic}
	If a function $f \in L^p(\mathbb{R})$, then $H \left (  (1- \chi_{(a,b)}) f \right )(x)$ is real analytic on any closed subinterval of $ (a,b)$.
\end{lemma}
\pf Take any $x \in [a_0, b_0] \subset (a,b)$, then
\begin{eqnarray*}
	H \left (  (1- \chi_{(a,b)}) f \right )(x) & = & \frac{1}{\pi} \int_{-\infty}^a \frac{f(y)}{x-y} \, dy + \frac{1}{\pi} \int_b^{\infty} \frac{f(y)}{x-y} \, dy . 
\end{eqnarray*}


\noindent
If $f \in L^p(\mathbb{R})$ for some $p>1$, then
\begin{eqnarray*}
	\left \vert	 \int_{-\infty}^a  \frac{\partial ^j}{\partial x^j}  \left ( \frac{f(y)}{x-y} \right ) \, dy  \right \vert & \leq & j! \int_{-\infty}^a \left \vert \frac{f(y)}{(x-y)^{j+1}}  \right \vert \, dy\\
	& \leq & j! \| f  \|_{L^p(\mathbb{R})}  \left ( \int_{-\infty}^a \left \vert \frac{1}{(x-y)^{(j+1)q}}  \right \vert \, dy \right )^\frac{1}{q} \hspace{.4in} ( \text{where} \,\, \frac{1}{p} + \frac{1}{q} =1  ) \\
	& = & j!\| f  \|_{L^p(\mathbb{R})}\dfrac{ (x-a)^{\frac{1-(j+1)q}{q}}}{\big (  (j+1)q-1\big )^{\frac{1}{q}}}  \\
	& \leq & j!\| f  \|_{L^p(\mathbb{R})}\left ( \frac{1}{a_0-a} \right )^{j+1}. 
\end{eqnarray*} 
If $p=1$, this is easily true.
Then by the Dominated Convergence Theorem for all $j \geq 1$,
\begin{equation*}
\left |	\frac{\partial ^j}{\partial x^j}  \left ( \int_{-\infty}^a \frac{f(y)}{x-y} \, dy \right ) \right | \,\, = \,\, \left | \int_{-\infty}^a  \frac{\partial ^j}{\partial x^j}   \left ( \frac{f(y)}{x-y} \right ) \, dy \right |  \,\, \leq \,\,   j!\| f  \|_{L^p(\mathbb{R})}\left ( \frac{1}{a_0-a} \right )^{j+1},
\end{equation*}

\noindent
thus $\displaystyle \int_{-\infty}^a \frac{f(y)}{x-y} \, dy$ is real analytic at $x$.  Similarly, we can prove that $\displaystyle \int_b^{\infty} \frac{f(y)}{x-y} \, dy $ is also real analytic at $x$, therefore $H \left (  (1- \chi_{[a,b]}) f \right ) $ is real analytic at $x$.

\stop

\vspace{.2in}

\noindent
\textit{Proof of Corollary \ref{cor:C_infinity}}:
If $f(x)$ is $C^{k+1}$ near a point, then in particular it is $C^{k, \alpha}$ for any $0<\alpha<1 $.  By Theorem \ref{thm:main} it follows that $Hf(x)$ is $C^{k, \alpha_0}$ for any $0<\alpha_0 <\alpha$ near the same point, therefore $Hf$ is $C^k$.  Furthermore, if $f(x)$ is $C^{\infty}$, then $Hf(x)$ is also $C^{\infty}$.  

\stop

\vspace{.2in}

\section{Proof of Theorem \ref{thm:real_analytic}}
\label{section:real_proof}

\vspace{.1in}
Now suppose $f$ is real analytic at a point $x_0$, we will prove that $Hf$ is also real analytic at $x_0$.
Since $f$ is real analytic at $x_0$, there is a $\delta_1 >0$ such that for any $y \in \left [ x_0-\delta_1 , x_0+\delta_1 \right ]$, 
$$f(y) = \sum_{l=0}^{\infty} \frac{f^{(l)}(x_0)}{l!}(y-x_0)^l, $$ where $$ |f^{(l)}(x_0)| \leq l!c^{l} $$ for some constant  $c$.
Let $$\delta_2 = \frac{1}{2c}$$ and $$\delta = \min \{\delta_1, \delta_2 \}.$$

\vspace{.05in}
\noindent
Denote $a=x_0 - \delta$ and $b=x_0+\delta$, again we write
$$Hf = H(\chi_{(a,b)}f) + H \left ( \left( 1-\chi_{(a,b)}\right )f \right ).$$
    Since $H \left ( \left( 1-\chi_{(a,b)}\right )f \right )$ is real analytic at $x_0$ by Lemma \ref{lemma:outside_analytic}, we 
only need to show  $H \left ( \chi_{(a,b)} f  \right )$ is also real analytic at $x_0$. By (\ref{eqn:H_chi_f})
$$ 	H \left ( \chi_{(a,b)} f \right )(x) = -\frac{1}{\pi}Tf(x)  + \frac{f(x)}{\pi}\ln \left \vert \frac{x-a}{x-b} \right \vert. $$
$f(x)$ being real analytic at $x_0$ implies $ \frac{f(x)}{\pi}\ln \left \vert \frac{x-a}{x-b} \right \vert$ being real analytic at $x_0$, so it remains to prove that $Tf(x)$ is real analytic at $x_0$.
By (\ref{eqn:Tf_kth_derivative})

\begin{eqnarray*}
	\left | \frac{d^k}{dx^k} Tf(x_0) \right | & = & \left | k!\int _a^b \frac{f(y)-P_k(x_0,y)}{(y-x_0)^{k+1}} dy \right | \\
	&  = & \left |  k!\int _a^b \sum_{l=k+1}^{\infty} \frac{f^{(l)}(x_0)}{l!}(y-x_0)^{l-k-1} dy \right | \\
	 & = & \left | k! \sum_{l=k+1}^{\infty} \dfrac{f^{(l)}(x_0)}{l!(l-k)} \left [ \delta^{l-k} - (-\delta)^{l-k} \right ] \right | \\
	 & \leq  & 2k! \sum_{l=k+1}^{\infty} \dfrac{c^{l}\delta^{l-k}}{l-k}.
\end{eqnarray*}

\noindent
Since $c\delta \leq \frac{1}{2}$, we have $$\sum_{l=k+1}^{\infty} \dfrac{c^{l}\delta^{l-k}}{l-k} <  \sum_{l=k+1}^{\infty} \left ( \frac{1}{2} \right )^l\delta^{-k} < \frac{1}{(2\delta)^k}.$$ 
Hence $$ 	\left| \frac{d^k}{dx^k} Tf(x_0) \right | \leq 2k!\left ( \frac{1}{2\delta}  \right )^k, $$ which implies $Tf$ is real analytic at $x_0$. 

\stop
\vspace{.2in}

\section{Examples for H\"older Sharpness}
\label{section:example}

\vspace{.1in}
In this section we present two examples where the functions do not lose H\"older regularity under the Hilbert transform.  The first example is the following.

\begin{prop}
	\label{prop:regularity_same}
	Let 	\begin{equation*}
		f(x) = 
		\begin{cases}
			\sqrt{ 1-x^2 } & \,\, \text{if} \,\, |x| \leq 1,  \\
			\noalign{\medskip}
			0	 & \,\, \text{if} \,\, |x|>1,
		\end{cases}
	\end{equation*}
	then 	
	\allowdisplaybreaks
	\begin{equation*}
		Hf(x) =
		\begin{cases}
			x & \,\, \text{if} \,\, |x| \leq 1,  \\
			\noalign{\medskip}
			x - \sgn(x)\sqrt{x^2-1}	 & \,\, \text{if} \,\,  |x|>1,
		\end{cases}
	\end{equation*}
	
	\noindent
	and both $f $ and $Hf$ are in $C^{0,\frac{1}{2}}\left ( \mathbb{R} \right ) \cap L^2(\mathbb{R}) \cap L^{\infty}(\mathbb{R}) . $
\end{prop}


\vspace{.1in}

To prove this example we need to use the following integral formula.  It was originally computed in \cite{Pandey} (P. 80 - 81) with an error in the case $x<-1$ where the minus sign was absent.

\begin{lemma}
		\label{lem:integral_minus_1_to_1}
\begin{equation}
	p.v. \,	\int _{-1}^1  \frac{1}{\sqrt{1-t^2} \, (x-t)}  \, dt = 
	\begin{cases}
		0 & \,\, \text{if} \,\, |x|<1,  \\
		\noalign{\medskip}
		\frac{\pi}{\sqrt{x^2-1}}  & \,\, \text{if} \,\, x >1, \\
		\noalign{\medskip}
		-\frac{\pi}{\sqrt{x^2-1}}  & \,\, \text{if} \,\, x < -1.
	\end{cases}
 \end{equation}

\end{lemma}

\vspace{.1in}

\pf
Letting $$t = \frac{1-y^2}{1+y^2}, \hspace{.2in} y>0,$$
then 
$$ \int _{-1}^1  \frac{1}{\sqrt{1-t^2} \, (x-t)}  \, dt =   2 \int_0^{\infty} \dfrac{1}{(1+x)y^2-(1-x)} \, dy.   $$

\vspace{.1in}
\noindent
If $|x|<1$, then $\frac{1+x}{1-x}>0$, thus
\allowdisplaybreaks
\begin{eqnarray*}
& & 2 \int_0^{\infty} \dfrac{1}{(1+x)y^2-(1-x)} \, dy \\
 & = & \frac{2}{1+x } \,\, p.v. \,\, \int_0^{\infty} \dfrac{1}{y^2-\frac{1-x}{1+x} } \, dy \\
& = &\frac{2}{1+x } \lim _{\epsilon \to 0}  \left ( \int_0^{\sqrt{\frac{1-x}{1+x}} - \epsilon} \dfrac{1}{y^2-\frac{1-x}{1+x} } \, dy + \int_{\sqrt{\frac{1-x}{1+x}} + \epsilon}^{\infty} \dfrac{1}{y^2-\frac{1-x}{1+x} } \, dy  \right ) \\
& = & \frac{2}{\sqrt{1-x^2} } \lim_{\epsilon \to 0} \left (  \ln \left \vert \frac{y- \sqrt{\frac{1-x}{1+x}}}{y+\sqrt{\frac{1-x}{1+x}}} \right \vert  \Bigg \vert_0^{\sqrt{\frac{1-x}{1+x}} - \epsilon} +   \ln \left \vert \frac{y- \sqrt{\frac{1-x}{1+x}}}{y+\sqrt{\frac{1-x}{1+x}}} \right \vert  \Bigg \vert_{\sqrt{\frac{1-x}{1+x}} + \epsilon}^{\infty}  \right ) \\
& = &\frac{2}{\sqrt{1-x^2} } \lim_{\epsilon \to 0} \left ( \ln \left \vert  \frac{\epsilon}{ 2\sqrt{\frac{1-x}{1+x}} -\epsilon  } \right  \vert  -  \ln \left \vert  \frac{\epsilon}{ 2\sqrt{\frac{1-x}{1+x}} + \epsilon  } \right  \vert \right ) \\
& = & \frac{2}{\sqrt{1-x^2} } \lim_{\epsilon \to 0} \ln \left \vert  \frac{ 2\sqrt{\frac{1-x}{1+x}} + \epsilon}{ 2\sqrt{\frac{1-x}{1+x}} -\epsilon  } \right  \vert   \\
& = & 0.
\end{eqnarray*}

\vspace{.1in}
\noindent
If $|x|>1$, then $ \frac{x-1}{x+1} >0 $, thus
 \allowdisplaybreaks
 \begin{eqnarray*}
 	 2 \int_0^{\infty} \dfrac{1}{(1+x)y^2-(1-x)} \, dy 
 	& = & \frac{2}{1+x }  \int_0^{\infty} \dfrac{1}{y^2+\frac{x-1}{x+1} } \, dy \\
 	& = & \frac{2}{1+x } \cdot \sqrt{\frac{x+1}{x-1}} \tan^{-1} \left ( \frac{y}{ \sqrt{\frac{x-1}{x+1}}  } \right ) \Bigg \vert_{0}^{\infty} \\
 	& = & 	\begin{cases}
 		\frac{\pi}{\sqrt{x^2-1}} & \,\, \text{if} \,\, x > 1,  \\
 		\noalign{\medskip}
 			-\frac{\pi}{\sqrt{x^2-1}} & \,\, \text{if} \,\, x < - 1, 
 	\end{cases}
 \end{eqnarray*}

\stop

\vspace{.2in}

\noindent
\textbf{Proof of Proposition \ref{prop:regularity_same}}: By definition $f \in C^{0, \frac{1}{2}} \left ( \mathbb{R} \right ) \cap L^2(\mathbb{R})\cap L^{\infty}(\mathbb{R}) $. Next we compute $Hf$.

\vspace{.1in}

\noindent
If $|x|<1$, then
\allowdisplaybreaks
\begin{eqnarray*}
			Hf(x) & = & p.v. \, \frac{1}{\pi } \int_{-\infty}^{\infty} \frac{f(t)}{x-t} \,\, dt \\
	& = & 	p.v. \, \frac{1}{\pi} \int _{-1}^{1} \frac{ \sqrt{ 1-t^2 } }{x-t} \, dt \\
	& = &  p.v. \, \frac{1}{\pi} \int _{-1}^{1} \frac{ 1-t^2  }{ \sqrt{ 1-t^2 } \, (x-t)} \, dt \\
	& = &  p.v. \, \frac{1}{\pi} \int _{-1}^{1} \frac{ 1  }{ \sqrt{ 1-t^2 } \, (x-t)} \, dt -   p.v. \, \frac{1}{\pi} \int _{-1}^{1} \frac{ t^2  }{ \sqrt{ 1-t^2 } \, (x-t)} \, dt . 
\end{eqnarray*}
By Lemma \ref{lem:integral_minus_1_to_1}, the first term $$p.v. \, \frac{1}{\pi} \int _{-1}^{1} \frac{ 1  }{ \sqrt{ 1-t^2 } \, (x-t)} \, dt =0,$$ and the second term can be computed as
\begin{eqnarray*}
& & p.v. \, \frac{1}{\pi} \int _{-1}^{1} \frac{ t^2  }{ \sqrt{ 1-t^2 } \, (x-t)} \, dt \\
& = & p.v. \, \frac{-1}{\pi} \int _{-1}^{1} \frac{ x^2 - t^2  }{ \sqrt{ 1-t^2 } \, (x-t)} \, dt \\
 & = & \lim_{\epsilon \to 0} \left ( -  \frac{1}{\pi} \int _{-1}^{x-\epsilon} \frac{x^2 - t^2  }{ \sqrt{ 1-t^2 } \, (x-t)} \, dt - \frac{1}{\pi} \int _{x+\epsilon}^{1} \frac{ x^2 - t^2  }{ \sqrt{ 1-t^2 } \, (x-t)} \, dt \right ) \\
  & = & \lim_{\epsilon \to 0} \left ( - \frac{1}{\pi} \int _{-1}^{x-\epsilon} \frac{x+t  }{ \sqrt{ 1-t^2 } } \, dt - \frac{1}{\pi} \int _{x+\epsilon}^{1} \frac{ x+t  }{ \sqrt{ 1-t^2 } } \, dt \right ) \\
& = &  \, - \frac{1}{\pi} \int _{-1}^{1} \frac{ x+t }{ \sqrt{ 1-t^2 } } \, dt \\
& = &  - \frac{1}{\pi} \int _{-1}^{1} \frac{ x }{ \sqrt{ 1-t^2 } } \, dt -  \frac{1}{\pi} \int _{-1}^{1} \frac{ t }{ \sqrt{ 1-t^2 } } \, dt \\
& = & - x.
\end{eqnarray*}
Therefore, $$ 	Hf (x) = x.$$

\vspace{.1in}
\noindent
If $|x|>1$, then 
\allowdisplaybreaks
\begin{eqnarray*}
	Hf (x) & = &   p.v. \, \frac{1}{\pi} \int_{-\infty}^{\infty} \frac{ f(t) }{x-t} \, dt\\
	& = & 	  \frac{1}{\pi} \int _{-1}^{1} \frac{ \sqrt{ 1-t^2 } }{x-t} \, dt \hspace{.9in}  (x-t \neq 0 \,\,\,\, \text{for all} \,\, -1 \leq t \leq 1) \\
	& = &    \frac{1}{\pi} \int _{-1}^{1} \frac{ 1-t^2  }{ \sqrt{ 1-t^2 } \, (x-t)} \, dt \\
	& = &   \frac{1}{\pi} \int _{-1}^{1} \frac{ 1  }{ \sqrt{ 1-t^2 } \, (x-t)} \, dt -  \frac{1}{\pi} \int _{-1}^{1} \frac{ t^2  }{ \sqrt{ 1-t^2 } \, (x-t)} \, dt . \\
\end{eqnarray*}

\noindent
By Lemma \ref{lem:integral_minus_1_to_1},  $$ \frac{1}{\pi} \int _{-1}^{1} \frac{ 1  }{ \sqrt{ 1-t^2 } \, (x-t)} \, dt = \frac{\sgn(x) }{\sqrt{x^2-1}} ,$$ and
\begin{eqnarray*}
	 \frac{1}{\pi} \int _{-1}^{1} \frac{ t^2  }{ \sqrt{ 1-t^2 } \, (x-t)} \, dt 
	& = & - \frac{1}{\pi} \int _{-1}^{1} \frac{ x^2 - t^2  }{ \sqrt{ 1-t^2 } \, (x-t)} \, dt +  \frac{1}{\pi} \int _{-1}^{1} \frac{ x^2 }{ \sqrt{ 1-t^2 } \, (x-t)} \, dt  \\
	& = & - \frac{1}{\pi} \int _{-1}^{1} \frac{ x + t  }{ \sqrt{ 1-t^2 } } \, dt + \frac{x^2 \sgn(x)}{\sqrt{x^2-1}} \\
	& = & - \frac{1}{\pi} \int _{-1}^{1} \frac{ x }{ \sqrt{ 1-t^2 } } \, dt -  \frac{1}{\pi} \int _{-1}^{1} \frac{ t }{ \sqrt{ 1-t^2 } } \, dt + \frac{x^2 \sgn(x)}{\sqrt{x^2-1}} \\
	 & =  & -\frac{x}{\pi} \sin^{-1} t \Big \vert_{t=-1}^{t=1} + 0 + \frac{x^2 \sgn(x)}{\sqrt{x^2-1}}\\
	& = & - x  + \frac{x^2 \sgn(x)}{\sqrt{x^2-1}}.
\end{eqnarray*}

\noindent
Therefore, 
\begin{eqnarray*} 
		Hf (x) & = & \frac{\sgn(x)}{\sqrt{x^2-1}} + x - \frac{x^2\sgn(x)}{\sqrt{x^2-1}} \\
		& = & x - \sgn(x)\sqrt{x^2-1}.
	\end{eqnarray*}

\vspace{.1in}
\noindent
Lastly, we need to compute the Hilbert transform of $f$ at $x=\pm 1$.  When $x=1$, 
\allowdisplaybreaks
\begin{eqnarray*}
		Hf (1) & = &   p.v. \, \frac{1}{\pi} \int_{-\infty}^{\infty} \frac{ f(t)  }{1-t} \, dt\\
	& = &   \frac{1}{\pi} \int _{-1}^{1} \frac{ \sqrt{ 1-t^2 } }{1-t} \, dt  \\
	& = &    \frac{1}{\pi} \int _{-1}^{1} \sqrt{\frac{1+t}{1-t}}\, dt. 
\end{eqnarray*}

\noindent
Let $s=\frac{1+t}{1-t}$, then $dt=\frac{2}{(s+1)^2}ds$, so
\allowdisplaybreaks
\begin{eqnarray*}
	 \frac{1}{\pi} \int _{-1}^{1} \sqrt{\frac{1+t}{1-t}}\, dt & = & \frac{1}{\pi} \int_0^{\infty} \frac{2\sqrt{s}}{(s+1)^2} ds\\
	& = &  \frac{1}{\pi} \int_1^{\infty} \frac{2\sqrt{s-1}}{s^2} ds\\
	& = & \frac{1}{\pi} \int_0^{\frac{\pi}{2}} 4\sin^2 \theta \, d\theta   \hspace{.6in} (\text{let} \hspace{.05in} s=\sec^2 \theta) \\
	& = & 1.
\end{eqnarray*}

\noindent
By similar calculations we can also prove 
$Hf (-1) = -1.  $
Therefore we have proved that
	\allowdisplaybreaks
\begin{equation*}
	Hf(x) =
	\begin{cases}
		x & \,\, \text{if} \,\, |x| \leq 1,  \\
		\noalign{\medskip}
		x - \sgn(x)\sqrt{x^2-1}	 & \,\, \text{if} \,\,  |x|>1.
	\end{cases}
\end{equation*}

\vspace{.05in}
\noindent
By straightforward but lengthy calculations it can be verified that $Hf(x) \in C^{0, \frac{1}{2}} \left ( \mathbb{R} \right ) \cap L^2 \left (\mathbb{R} \right ) \cap L^{\infty}(\mathbb{R}) $.  We will provide the details in Appendix B.

\stop

\vspace{.2in}

To construct our second example, we will start with the simple function $x^{\alpha}$ and compute its finite Hilbert transform first.

\begin{lemma}
	\label{lemma:x^alpha_finite_Hilbert}
	For any $0<x<1$, the finite Hilbert transform of the function $x^{\alpha}$ on $(0,1)$ is 
\begin{equation*}
	p.v. \, \frac{1}{\pi}\int_0^1 \frac{y^{\alpha}}{x-y} dy =  - \frac{(1-x)^{\alpha}}{\pi \alpha} - \frac{1}{\pi}  x^{\alpha}g(x),
\end{equation*}
where $g(x)$ is a continuous function on $[0,1]$ and $g(0) \neq 0$. 
\end{lemma}

\pf
$$ \int _0 ^1 \dfrac{y^{\alpha} - x^{\alpha}  }{y-x} \, dy \,\, = \int _0 ^x  \dfrac{y^{\alpha} - x^{\alpha}  }{y-x} \, dy + \int _x ^1  \dfrac{y^{\alpha} - x^{\alpha}  }{y-x} \, dy\,\, := \,\,  I_1 + I_2. $$

\vspace{.1in}
\noindent
We first compute $I_1$.  Let $z=\frac{y}{x}$, then $dy = x dz$ and $0<z<1$ when $0<y<x$. It follows that
\begin{equation}
	\label{eqn:int_0_to_x}
	I_1 \,\, = \,\, 	\int _0 ^x  \dfrac{y^{\alpha} - x^{\alpha}  }{y-x} \, dy \,\, = \,\,  x^{\alpha} 	\int _0 ^1  \dfrac{  z^{\alpha}-  1 }{z-1} \, dz.
\end{equation}	

\noindent
Since
$$ \lim _{z \to 1} \dfrac{ z^{\alpha}- 1 }{z-1} = \alpha, $$
\noindent
the function $  \dfrac{ z^{\alpha}- 1 }{z-1} $ has a removable discontinuity at $1$, so $$ \displaystyle 	\int _0 ^1  \dfrac{  z^{\alpha}-  1 }{z-1} \, dz= C_1 \,\,\,\, \text{for some constant}  \,\, C_1>0.$$  Hence, we can write 
\begin{equation}
	\label{eqn:estimate_I_1}
	I_1 = C_1x^{\alpha}
\end{equation}    

\noindent
Next we compute $I_2$.  Again let $z=\frac{y}{x}$, then $dy = x dz$ and $1< z < \frac{1}{x} $ when $x<y<1$. Similar to the previous case we can write
\begin{equation}
	\label{eqn: int_x_to_1}
	I_2 \,\, = \,\,	\int _x ^1  \dfrac{y^{\alpha} - x^{\alpha}  }{y-x} \, dy \,\, = \,\,  x^{\alpha} 	\int _1 ^{\frac{1}{x}}  \dfrac{  z^{\alpha}-  1 }{z-1} \, dz.
\end{equation}	  

\noindent
As $x \to 0$, $\frac{1}{x} \to \infty$, so we only consider small $x$ such that $\frac{1}{x} >2$.  Since $\dfrac{  z^{\alpha}-  1 }{z-1}$ has a removable discontinuity at $1$,  $$ \displaystyle \int _1 ^2  \dfrac{  z^{\alpha}-  1 }{z-1} \, dz= C_2 \,\,\,\, \text{for some constant}  \,\, C_2>0.$$  Hence, we can write 
\begin{equation}
	\label{eqn:estimate_integral_1_to_2}
	x^{\alpha} \int _1 ^2   \dfrac{  z^{\alpha}-  1 }{z-1}\, dz  = C_2x^{\alpha}.
\end{equation}    

\vspace{.1in}
\noindent
It now remains to compute the integral from 2 to $\frac{1}{x}$.

\begin{equation}
	\label{eqn:split_integrals}
	\int _2 ^{\frac{1}{x}} \dfrac{  z^{\alpha}-  1 }{z-1}  \, dz \,\, = \,\, \int _2 ^{\frac{1}{x}}  \dfrac{  z^{\alpha} }{z-1} \, dz + \ln \left \vert \frac{x}{1-x} \right \vert.
\end{equation}

\vspace{.1in}
\noindent
To evaluate the first term, 
let $$ t = \frac{1}{z-1},$$ so $0<t<1$ when $z > 2$, and $$  z=1+\frac{1}{t} \,\, \,\, \text{with} \,\,\,\, dz=-\frac{1}{t^2}dt. $$
Then 
\begin{equation}
	\label{eqn:int_y_to_t}
	\int _2 ^{\frac{1}{x}}   \dfrac{ z^{\alpha} }{z-1} \, dz =  \int_{\frac{x}{1-x}}^1 \frac{\left ( 1+t \right )^{\alpha}}{t^{1+\alpha}} \, dt.
\end{equation} 

\vspace{.1in}
\noindent
Since $0<t<1$, 
$$\left ( 1+ t \right )^{\alpha} = 1+ a_1 t + a_2t^2 + a_3t^3\cdots, $$ where $a_1=\alpha  $ and in general $$a_i = \dfrac{\alpha(\alpha -1) \cdots (\alpha - i+1)} {i!} \,\,\,\, \text{for} \,\,\,\, i \geq 1.$$

\vspace{.1in}
\noindent
It follows that 
\begin{equation}
	\label{eqn:term_by_term_t_integral}
	\int_{\frac{x}{1-x}}^1 \frac{\left ( 1+t \right )^{\alpha}}{t^{1+\alpha}} \, dt  =   \int_{\frac{x}{1-x}}^1  \left (t^{-1-\alpha} + a_1 t^{-\alpha} + a_2t^{1-\alpha} + a_3t^{2-\alpha} + \cdots \right ) \, dt.
\end{equation}

\vspace{.1in}
\noindent
Since the binomial series converges absolutely, we can integrate term-by-term to obtain
\begin{eqnarray*}
	& & \int  \left (t^{-1-\alpha} + a_1 t^{-\alpha} + a_2t^{1-\alpha} + a_3t^{2-\alpha} + \cdots \right ) \, dt \\
	& = &  \int  t^{-1-\alpha} \, dt + \int a_1 t^{-\alpha} \, dt + \int a_2 t^{1-\alpha} \, dt + \cdots \\
	& = & \frac{t^{-\alpha}}{-\alpha} + \frac{a_1t^{1-\alpha}}{1-\alpha} + \frac{a_2 t^{2-\alpha}}{2-\alpha} + \cdots \\
	& = & \frac{t^{-\alpha}}{-\alpha} + \sum _{k=1}^{\infty}  \frac{a_k }{k-\alpha}t^{k-\alpha}. 
\end{eqnarray*}

\vspace{.1in}
\noindent
Note that
$$\lim_{k \to \infty} \left | \dfrac{\,\,\,\, \frac{a_{k+1}}{k+1-\alpha}t^{k+1-\alpha} \,\,\,\, }{\frac{a_k}{k-\alpha}t^{k-\alpha}}  \right | = \lim_{k \to \infty} \left | \frac{\alpha -k}{k+1} \cdot \frac{k-\alpha}{k+1-\alpha} \cdot t \right | = t , $$ 

\noindent
thus $\displaystyle \sum _{k=1}^{\infty}  \frac{a_k }{k-\alpha}t^{k-\alpha}$ converges absolutely when $|t|<1$.  When $x$ is sufficiently close to $0$, we have $0<\frac{x}{1-x} <1,$ therefore $$  \sum _{k=1}^{\infty}  \frac{a_k }{k-\alpha} -  \sum _{k=1}^{\infty}  \frac{a_k }{k-\alpha}\left ( \frac{x}{1-x} \right )^{k-\alpha} = \sum _{k=1}^{\infty}  \frac{a_k }{k-\alpha} \left [ 1- \left ( \frac{x}{1-x} \right )^{k-\alpha}\right ] . $$

\vspace{.1in}
\noindent
Thus we can write 
\begin{eqnarray*}
	& & 	\int_{\frac{x}{1-x}}^1  \left (t^{-1-\alpha} + a_1 t^{-\alpha} + a_2t^{1-\alpha} + a_3t^{2-\alpha} + \cdots \right ) \, dt \\
	& = &  \left ( \frac{t^{-\alpha}}{-\alpha} + \sum _{k=1}^{\infty}  \frac{a_k }{k-\alpha}t^{k-\alpha} \right ) \bigg \vert _{t=\frac{x}{1-x}}^{t=1} \\
	& = &  \frac{(1-x)^{\alpha}x^{-\alpha} - 1}{\alpha} + \sum _{k=1}^{\infty}  \frac{a_k }{k-\alpha} \left [ 1- \left ( \frac{x}{1-x} \right )^{k-\alpha}\right ] .
\end{eqnarray*}

\vspace{.1in}
\noindent
Together with (\ref{eqn:term_by_term_t_integral}) this implies
\begin{eqnarray}
	\label{eqn:int_2nd_term}
	&  &  \int_{\frac{x}{1-x}}^1 \frac{\left ( 1+t \right )^{1+\alpha}}{t^{1+\alpha}} \, dt \\
	& = & \frac{(1-x)^{\alpha}x^{-\alpha} - 1}{\alpha} + \frac{a_1}{1-\alpha} \left ( 1 - \dfrac{x^{1-\alpha}}{(1-x)^{1-\alpha}}\right ) + \frac{a_2}{2-\alpha } \left (  1 - \dfrac{x^{2-\alpha}}{(1-x)^{2-\alpha}}\  \right ) + \cdots . \nonumber
\end{eqnarray}

\vspace{.1in}
\noindent
By (\ref{eqn:split_integrals}), (\ref{eqn:int_y_to_t}) and (\ref{eqn:int_2nd_term}), we have
\begin{eqnarray}
	\label{eqn:compute_int_2_to_1/x_with_x_alpha}
	& & x^{\alpha}\int _2 ^{\frac{1}{x}}   \dfrac{ z^{\alpha}- 1 }{z-1} \, dz \\
	& = &  x^{\alpha} \ln \left \vert \frac{x}{1-x} \right \vert + \frac{(1-x)^{\alpha}- x^{\alpha}}{\alpha} \nonumber \\
	& + &  \frac{a_1}{1-\alpha} \left ( x^{\alpha} - \dfrac{x}{(1-x)^{1-\alpha}}\right ) + \frac{a_2}{2-\alpha } \left (  x^{\alpha} - \dfrac{x^{2}}{(1-x)^{2-\alpha}}\  \right ) + \cdots . \nonumber
\end{eqnarray}

\vspace{.1in}
\noindent
Combined with (\ref{eqn:estimate_I_1}), (\ref{eqn: int_x_to_1}), and (\ref{eqn:estimate_integral_1_to_2}), this implies that
\allowdisplaybreaks
\begin{eqnarray}
	\label{eqn:F(x)_estimate}
	\int _0 ^1 \dfrac{y^{\alpha} - x^{\alpha}  }{y-x} \, dy & = & \left ( C_1+C_2 \right ) x^{\alpha} +  x^{\alpha}\ln \left \vert \frac{x}{1-x} \right \vert +  \frac{(1-x)^{\alpha}}{\alpha}  -  \frac{x^{\alpha}}{\alpha} \nonumber \\
	&  & + \frac{a_1}{1-\alpha} \left ( x^{\alpha} - \dfrac{x}{(1-x)^{1-\alpha}}\right ) + \frac{a_2}{2-\alpha } \left (  x^{\alpha} - \dfrac{x^{2}}{(1-x)^{2-\alpha}}\  \right ) + \cdots \nonumber \\
	& = &  x^{\alpha}\ln \left \vert \frac{x}{1-x} \right \vert +  \frac{(1-x)^{\alpha}}{\alpha} + x^{\alpha}g(x),
\end{eqnarray}
where $g(x)$ is a continuous function on $[0,1]$, and it can be verified that $g(0) \neq 0$.    
By direct calculation $$ p.v.\,\frac{1}{\pi} \int_{0}^{1} \frac{ x^{\alpha}}{x-y} dy = \frac{1}{\pi} x^{\alpha} \ln \left \vert \frac{x}{1-x} \right \vert .$$ 

\noindent
Therefore the finite Hilbert transform of $x^{\alpha}$ is
\begin{eqnarray*}
	p.v. \, \frac{1}{\pi}\int_0^1 \frac{y^{\alpha}}{x-y} dy & = & - \frac{1}{\pi}	\int _0 ^1 \dfrac{y^{\alpha} - x^{\alpha}  }{y-x} \, dy + p.v. \, \frac{1}{\pi} \int _0 ^1 \dfrac{x^{\alpha}  }{x-y} \, dy \\	
	& = & - \frac{(1-x)^{\alpha}}{\pi \alpha} - \frac{1}{\pi}  x^{\alpha}g(x).
\end{eqnarray*}

\vspace{.1in}
\stop

\vspace{.1in}
\noindent
\textbf{Remark}: The logarithmic term in Equation (\ref{eqn:F(x)_estimate}) shows that the function $x^{\alpha}$ loses regularity at $x=0$ under the operation $T$, therefore Theorem \ref{thm:holder} (and similarly Theorem \ref{thm:T_k_Holder}) cannot be improved to $\alpha_0 = \alpha$.  On the other hand, Lemma \ref{lemma:x^alpha_finite_Hilbert} shows that there is no loss of regularity under the finite Hilbert transform of $x^{\alpha}$, therefore we cannot yet rule out the possibility of Theorem \ref{thm:main} being true for $\alpha_0 = \alpha$.  However, if it is true, a more powerful tool than the operator $T$ must be employed to prove it.

\vspace{.2in}

Next, we will shift $x^\alpha$ to an arbitrary point $x_0 \in (0,1)$ and show that its finite Hilbert transform is H\"older continuous with exponent $\alpha$ near $x_0$.  Combined with a cut-off function, it provides an example of a $C^{0, \alpha}$ function whose Hilbert transform is still locally $C^{0, \alpha}$.

\begin{prop}
	\label{prop:x_alpha}
	Given any $0 < \alpha <1$ and $x_0 \in (0,1)$, define
	\begin{equation*}
		\label{defn:f_shifted}
		f(x) = 
		\begin{cases}
			0 & \,\, \text{if} \,\,\,\,   x \leq x_0,  \\
			\noalign{\medskip}
			(x -x_0)^{\alpha}  & \,\, \text{if} \,\,\,\,x > x_0.  
		\end{cases}
	\end{equation*}
	Then the finite Hilbert transform of $f$ on $(0,1)$ is
	\begin{eqnarray*}
		p.v. \, \frac{1}{\pi} \int_0^1 \frac{f(y)}{x-y} \, dy 
		& = & -  \frac{(1-x)^{ \alpha}}{\pi \alpha} -    \frac{(x-x_0)^{\alpha}}{\pi} g\left (  \frac{x-x_0}{1-x_0} \right ),
	\end{eqnarray*}
	where $g(x)$ is a continuous function on $[0,1]$ and $g(0) \neq 0$.  Therefore the finite Hilbert transform of $f$ is H\"older continuous with exponent $\alpha$ in a neighborhood of $x_0$.

\end{prop}

\vspace{.1in}
\pf
By definition 
\begin{eqnarray*}
	\int_0^1 \frac{f(y) - f(x)}{y-x} \, dy 
	& = & - \int_0^{x_0} \frac{    \left ( x-x_0 \right )^{\alpha}}{y-x} \, dy + \int_{x_0}^1 \frac{ \left ( y-x_0 \right )^{\alpha} -  \left ( x-x_0 \right )^{\alpha}}{y-x} \, dy  .
\end{eqnarray*}

\noindent
The first term can be computed directly as
\begin{equation}
	\label{eqn:interior_1st_term}
	- \int_0^{x_0} \frac{ \left ( x-x_0 \right )^{\alpha}}{y-x} \, dy = \left ( x-x_0 \right )^{\alpha}  \ln \left \vert \frac{x}{x-x_0}  \right \vert.
\end{equation}

\noindent
To compute the second term, we make a change of variables by letting $$\hat{y} = \frac{y-x_0}{1-x_0}  \hspace{.2in} \text{and} \hspace{.2in}  \hat{x} = \frac{x-x_0}{1-x_0}.$$

\noindent
Then 
\begin{eqnarray*}
	&  &  \int_{x_0}^1 \frac{ \left ( y-x_0 \right )^{\alpha} -  \left ( x-x_0 \right )^{\alpha} }{y-x} \, dy \\
	& = & \left ( 1- x_0 \right )^{\alpha} \int_0 ^1 	\frac{\hat{y}^{\alpha} - \hat{x}^{\alpha} }{ \hat{y} - \hat{x} } \, d\hat{y} \\
	& = & \left ( 1- x_0 \right )^{\alpha} \left (  \hat{x}^{\alpha}\ln \left \vert \frac{\hat{x}}{1-\hat{x}} \right \vert +  \frac{(1-\hat{x})^{\alpha}}{\alpha} + \hat{x}^{\alpha}g(\hat{x}) \right ),
\end{eqnarray*}

\vspace{.1in}
\noindent
where in the last step we have used (\ref{eqn:F(x)_estimate}).
Written in the $x$ and $y$ variables the equation is
\begin{eqnarray}
	\label{eqn:interior_2nd_term}
	& &	\int_{x_0}^1 \frac{ \left ( y-x_0 \right )^{\alpha} -  \left ( x-x_0 \right )^{\alpha} }{y-x} \, dy  \\
	& = &  (x-x_0)^{\alpha} \ln \left \vert \frac{x-x_0}{1-x} \right \vert +  \frac{(1-x)^{\alpha}}{\alpha} +   (x-x_0)^{\alpha}g\left (  \frac{x-x_0}{1-x_0} \right ). \nonumber
\end{eqnarray}

\vspace{.1in}
\noindent
Combining (\ref{eqn:interior_1st_term}) and (\ref{eqn:interior_2nd_term}) we have 
\begin{equation*}
	\int_0^1 \frac{f(y) - f(x)}{y-x} \, dy  = (x-x_0)^{\alpha}  \ln \left \vert \frac{x}{1-x} \right \vert +\frac{(1-x)^{\alpha}}{\alpha} +   (x-x_0)^{\alpha}g\left (  \frac{x-x_0}{1-x_0} \right ).
\end{equation*}

\vspace{.1in}
\noindent
Thus the finite Hilbert tranform of $f$ on $(0,1)$ is
\begin{eqnarray*}
	&  & 	\int_0^1 \frac{f(y)}{x-y} \, dy \\
	& = &  - \int_0^1 \frac{f(y) - f(x)}{y-x} \, dy + \int_0^1 \frac{ f(x)}{x-y} \, dy \\
	& = & -(x-x_0)^{\alpha} \ln \left \vert \frac{x}{1-x} \right \vert - \frac{(1-x)^{\alpha}}{\alpha} -   (x-x_0)^{\alpha}g\left (  \frac{x-x_0}{1-x_0} \right ) + (x-x_0)^{\alpha} \ln \left \vert \frac{x}{1-x} \right \vert \\
	& = &   -  \frac{(1-x)^{ \alpha}}{ \alpha} -    (x-x_0)^{\alpha} g\left (  \frac{x-x_0}{1-x_0} \right ),
\end{eqnarray*}

\vspace{.1in}
\noindent
 and it is straightforward to verify that this is H\"older continuous with exponent $\alpha$ in a neighborhood of $x_0$.

 \stop

\vspace{.4in}

\section{Appendix A}
 \label{sec:appendix}

\vspace{.2in}

\noindent
In this section we give a proof of \textbf{Lemma \ref{lem:deriv_convergence}}, which is based on the hint for an exercise in \cite{Tao}.

\vspace{.1in}
\noindent
\textbf{\textbf{Lemma 2.6}}: {\textit{	Let $\{ f_n (x) \}$ be a sequence of differentiable functions from $[a,b] $ to $\mathbb{R}$. Suppose
		$f_n \to f$ uniformly on $[a,b]$ and
		$f'_n \to g$ uniformly on $[a,b]$, where $f,g: [a,b] \to \mathbb{R}$.
		Then $f'(x)=g(x)$ for every $x \in [a,b]$.}

\vspace{.2in}

\pf
We will prove $f'(x)=g(x)$ for every $x \in [a,b]$ by the definition of derivatives.  Since $\{f'_n\} \to g$ uniformly, for any $\epsilon >0$, there is $N_1 \in \mathbb{N}$ such that 
$$ \left | f'_n(x) - g(x) \right | < \epsilon$$ for all $n \geq N_1$ and $x \in [a,b]$.  For any $h \neq 0$ and $x, x+h \in [a,b]$,
\begin{eqnarray}
	& &	\left |  \dfrac{f_m(x+h)-f_m(x)}{h} - \dfrac{f_n(x+h)-f_n(x)}{h}   \right | \nonumber\\
	& = & \left |   \dfrac{(f_m-f_n)(x+h) - (f_m-f_n)(x)  }{h} \right | \nonumber \\
	& = & \left | (f'_m - f'_n) (\zeta) \right |,  \nonumber
\end{eqnarray} 
where $\zeta$ is a number between $x$ and $x+h$.  Since $f'_n \to g$ uniformly, there is $N_2 \in \mathbb{N}$ independent of $\zeta$ such that $$\left| f'_m (\zeta) - f'_n (\zeta) \right | < \epsilon$$ for all $m,n \geq N_2$. Thus
\begin{equation*}
	\left |  \dfrac{f_m(x+h)-f_m(x)}{h} - \dfrac{f_n(x+h)-f_n(x)}{h}   \right | < \epsilon
\end{equation*} 
for all $m, n \geq N_2$. Letting $m \to \infty$, this implies
\begin{equation*}
	\left |  \dfrac{f(x+h)-f(x)}{h} - \dfrac{f_n(x+h)-f_n(x)}{h}   \right | \leq \epsilon
\end{equation*} 
for all $n \geq N_2$. Now let $N=\max\{N_1, N_2\}$, then 
$$ \left | \dfrac{f(x+h)-f(x)}{h} - \dfrac{ f_N(x+h)-f_N(x)}{h} \right | \leq \epsilon$$ and 
$$  \left |  f'_N (x) - g(x) \right | < \epsilon $$
for all $x, x+h \in [a,b]$.
\vspace{.1in}

\noindent
For any fixed $x \in [a,b]$, by the definition of $f'_N(x)$, there is $\delta >0$ such that $$ \left | \dfrac{f_N(x+h)-f_N(x)}{h} - f'_N(x) \right | < \epsilon$$ for all $0<|h| < \delta$.  Then it follows that
\begin{eqnarray*}
	\left | \dfrac{ f(x+h)-f(x)}{h} -g(x) \right | 
	& \leq & \left | \dfrac{ f(x+h)-f(x)}{h} - \dfrac{f_N(x+h)-f_N(x)}{h} \right | \\
	& + & \left | \dfrac{f_N(x+h)-f_N(x)}{h} - f'_N(x) \right | + \left |  f'_N (x) - g(x) \right | \\
	&< & 3 \epsilon
\end{eqnarray*}
for all $0< |h| < \delta$, therefore $f'(x)=g(x)$.

\stop

\vspace{.2in}

\section{Appendix B}
\label{sec:Holder-example}

\vspace{.1in}
In the proof of Proposition \ref{prop:regularity_same} we showed that if \begin{equation*}
		f(x) = 
		\begin{cases}
			\sqrt{ 1-x^2 } & \,\, \text{if} \,\, |x| \leq 1,  \\
			\noalign{\medskip}
			0	 & \,\, \text{if} \,\, |x|>1,
		\end{cases}
	\end{equation*}
	then 	
	\allowdisplaybreaks
	\begin{equation*}
		Hf(x) =
		\begin{cases}
			x & \,\, \text{if} \,\, |x| \leq 1,  \\
			\noalign{\medskip}
			x - \sgn(x)\sqrt{x^2-1}	 & \,\, \text{if} \,\,  |x|>1.
		\end{cases}
	\end{equation*}

\noindent
In this appendix we will show that $Hf(x) \in C^{0,\frac{1}{2}}\left ( \mathbb{R} \right ) \cap L^2(\mathbb{R}) \cap L^{\infty}(\mathbb{R})$.
Since 
$$ \lim_{x \to \infty} Hf(x) = \lim_{x \to \infty} \left ( x - \sqrt{x^2-1} \right ) =  \lim_{x \to \infty} \frac{1}{x + \sqrt{x^2-1}} =0 $$
and
$$ \lim_{x \to -\infty} Hf(x) = \lim_{x \to -\infty} \left ( x + \sqrt{x^2-1} \right ) =  \lim_{x \to -\infty} \frac{1}{x - \sqrt{x^2-1}} =0, $$ 

\noindent
we know $Hf (x) \in L^{\infty}(\mathbb{R})$.
Since $$ \int_1^{\infty} \left ( x - \sqrt{x^2-1} \right )^2 \, dx = \int_1^{\infty} \frac{1}{\left (x + \sqrt{x^2-1} \right )^2 } \, dx <  \int_1^{\infty} \frac{1}{x^2}  <\infty $$ and
$$ \int_{-\infty}^{-1} \left ( x + \sqrt{x^2-1} \right )^2 \, dx = \int_{-\infty}^{-1} \frac{1}{\left (x - \sqrt{x^2-1} \right )^2 } \, dx <  \int_{-\infty}^{-1} \frac{1}{x^2}  <\infty, $$ we know $Hf(x) \in L^2(\mathbb{R})$.  Hence it remains to show $Hf(x) \in C^{0,\frac{1}{2}}\left ( \mathbb{R} \right )$.  Since $Hf$ is  bounded, we only need to estimate $ \frac{ \left | Hf(x_1) - Hf(x_2)  \right | }{\left | x_1 - x_2 \right |^{\frac{1}{2}}} $ when $| x_1 - x_2 | <1$.  There are five cases to consider.

\vspace{.1in}
\noindent
Case 1: $| x_1 - x_2 | <1$ and $x_1>x_2 \geq 1$.  We have
	\begin{eqnarray}
		\label{example_first_case}
		\frac{ \left | Hf(x_1) - Hf(x_2)  \right | }{\left | x_1 - x_2 \right |^{\frac{1}{2}}} & = &  	\frac{ \left | \left ( x_1 - \sqrt{x_1^2-1} \right ) - \left ( x_2 - \sqrt{x_2^2-1} \right )  \right | }{ \sqrt{x_1-x_2}} \nonumber \\
		& = & \left | \sqrt{x_1-x_2} - 	\frac{  \sqrt{x_1^2-1} - \sqrt{x_2^2-1} }{ \sqrt{x_1-x_2}} \right | . 
	\end{eqnarray}

	We will consider two separate subcases:  $1 \leq x_2 \leq 2$ or $x_2>2$.
	
	\begin{itemize}
		\item  Subcase 1: $1 \leq x_2 \leq 2$.  
		Let $t=x_1-x_2$, then 
		$x_1=x_2+t,$ $0<t<1$, and 
		\begin{eqnarray*}
		& & 	\left | \sqrt{x_1-x_2} - 	\frac{  \sqrt{x_1^2-1} - \sqrt{x_2^2-1} }{ \sqrt{x_1-x_2}} \right | \\
		& < & \sqrt{x_1-x_2} + \frac{  \sqrt{x_1^2-1} - \sqrt{x_2^2-1} }{ \sqrt{x_1-x_2}} \\
			& = & \sqrt{t} + \dfrac{\sqrt{x_2^2+2x_2t+t^2-1}-\sqrt{x_2^2-1}}{\sqrt{t}} \\
			& < & \sqrt{t} + \dfrac{ \left ( \sqrt{x_2^2-1} + \sqrt{2x_2t+t^2}  \right )-\sqrt{x_2^2-1}}{\sqrt{t}} \\
			& = & \sqrt{t} + \sqrt{2x_2+t} \\
			& < & 1+\sqrt{5 }.
		\end{eqnarray*}

			\item Subcase 2: $x_2 > 2$.  Then 
			\begin{eqnarray*}
& &	\left | \sqrt{x_1-x_2} - \dfrac{\sqrt{x_1^2-1}-\sqrt{x_2^2-1}}{\sqrt{x_1-x_2} } \right | \\ 
	& = & \sqrt{x_1-x_2} \left | 1- \dfrac{x_1+x_2 }{ \sqrt{x_1^2-1}+ \sqrt{x_2^2-1} }  \right | \\
	& < & \left | 1- \dfrac{x_1+x_2 }{ \sqrt{x_1^2-1}+ \sqrt{x_2^2-1} }  \right | .
\end{eqnarray*}

\noindent
The function $g(x_1, x_2)=\left | 1- \dfrac{x_1+x_2 }{ \sqrt{x_1^2-1}+ \sqrt{x_2^2-1} }  \right |$ is continuous for $x_1>x_2>2$, and $g(x_1, x_2) \to 0$ when $x_1, x_2 \to \infty$, thus $g$ is globally bounded.

\vspace{.1in}		
	\end{itemize}

In both subcases by (\ref{example_first_case}) it follows that  $ \frac{ \left | Hf(x_1) - Hf(x_2)  \right | }{\left | x_1 - x_2 \right |^{\frac{1}{2}}} $ is bounded for $x_1>x_2 \geq 1$.

\vspace{.1in}
\noindent
Case 2: $| x_1 - x_2 | <1$ and $x_2<x_1 \leq -1$.  We have
\allowdisplaybreaks 
\begin{eqnarray*}
	\frac{ \left | Hf(x_1) - Hf(x_2)  \right | }{\left | x_1 - x_2 \right |^{\frac{1}{2}}} & = & 	\frac{ \left | \left ( x_1 + \sqrt{x_1^2-1} \right ) - \left ( x_2 + \sqrt{x_2^2-1} \right )  \right | }{ \sqrt{x_1-x_2}} \\
	& = & \left | \sqrt{x_1-x_2} + \dfrac{\sqrt{x_1^2-1}-\sqrt{x_2^2-1}}{\sqrt{x_1-x_2} } \right |  .
\end{eqnarray*}

Denote $y_1=-x_2$ and $y_2=-x_1$.  Then $y_1>y_2 \geq 1$ and 
$$  \left | \sqrt{x_1-x_2} + \dfrac{\sqrt{x_1^2-1}-\sqrt{x_2^2-1}}{\sqrt{x_1-x_2} } \right |  = \left | \sqrt{y_1-y_2} - \dfrac{\sqrt{y_1^2-1}-\sqrt{y_2^2-1}}{\sqrt{y_1-y_2} } \right |,  $$

which is bounded by Case 1 above.
Therefore,  $ \frac{ \left | Hf(x_1) - Hf(x_2)  \right | }{\left | x_1 - x_2 \right |^{\frac{1}{2}}} $ is bounded for $x_2<x_1 \leq -1$.

\vspace{.1in}

\noindent
Case 3: $| x_1 - x_2 | <1$ and $-1 < x_2<x_1<-1$. We have
$$ 	\frac{ \left | Hf(x_1) - Hf(x_2)  \right | }{\left | x_1 - x_2 \right |^{\frac{1}{2}}} = \sqrt{x_1-x_2}<1. $$

\vspace{.1in}
\noindent
Case 4: $| x_1 - x_2 | <1$,  $-1 \leq  x_2 \leq 1 $, and $x_1 >1$. We have  $$ 	\frac{ \left | Hf(x_1) - Hf(x_2)  \right | }{\left | x_1 - x_2 \right |^{\frac{1}{2}}}  =  	\frac{ \left | \left ( x_1 - \sqrt{x_1^2-1} \right ) - x_2  \right | }{ \sqrt{x_1-x_2}} = 	\frac{ \left |  x_1 - x_2 - \sqrt{x_1^2-1}   \right | }{ \sqrt{x_1-x_2}} .$$

\vspace{.05in}
There are two subcases.

\vspace{.1in}

\begin{itemize}
	\item  Subcase 1: $x_1 - x_2 \geq \sqrt{x_1^2-1}$.  Then
	$$ 		\frac{ \left |  x_1 - x_2 - \sqrt{x_1^2-1}   \right | }{ \sqrt{x_1-x_2}} = 	\frac{  x_1 - x_2 - \sqrt{x_1^2-1}   }{ \sqrt{x_1-x_2}} < \frac{  x_1 - x_2   }{ \sqrt{x_1-x_2}} < \sqrt{x_1-x_2} <1. $$
	
	\vspace{.1in}
	
	\item 
Subcase 2: $0< x_1 - x_2 < \sqrt{x_1^2-1}$. Then 

\begin{equation}
	\label{eqn:example_holder_case 2}
		\frac{ \left |  x_1 - x_2 - \sqrt{x_1^2-1}   \right | }{ \sqrt{x_1-x_2}}  =  	\frac{   \sqrt{x_1^2-1} - ( x_1 - x_2)  }{ \sqrt{x_1-x_2}} \\
		<	\frac{   \sqrt{x_1^2-1}  }{ \sqrt{x_1-x_2}}. 
\end{equation}
From $| x_1 - x_2 | <1$,  $-1 \leq  x_2 \leq 1 $, and $x_1>1$ we know $1<x_1<2$.  Thus $(x_1-1)(x_1-3)<0$, i.e. $x_1^2-4x_1+3 <0$, which implies
$$ 0< x_1^2-1 < 4(x_1-1) ,$$ and consequently
$$ \sqrt{x_1^2-1} < 2 \sqrt{x_1-1}.  $$
Then  (\ref{eqn:example_holder_case 2}) implies 
$$	\frac{ \left |  x_1 - x_2 - \sqrt{x_1^2-1}   \right | }{ \sqrt{x_1-x_2}} < 	\frac{   \sqrt{x_1^2-1}  }{ \sqrt{x_1-x_2}} < 2 \frac{   \sqrt{x_1-1}  }{ \sqrt{x_1-x_2}} <  2.$$ 

\end{itemize}

Thus in both subcases $ 	\frac{ \left | Hf(x_1) - Hf(x_2)  \right | }{\left | x_1 - x_2 \right |^{\frac{1}{2}}} $ is bounded.

\vspace{.1in}

\noindent
Case 5: $| x_1 - x_2 | <1$,  $-1 \leq  x_1 \leq 1 $, and $x_2 < -1$.  We have $$	\frac{ \left | Hf(x_1) - Hf(x_2)  \right | }{\left | x_1 - x_2 \right |^{\frac{1}{2}}}  =  	\frac{ \left |x_1 - \left ( x_2 + \sqrt{x_2^2-1} \right )   \right | }{ \sqrt{x_1-x_2}} = 	\frac{ \left |  x_1 - x_2 - \sqrt{x_2^2-1}   \right | }{ \sqrt{x_1-x_2}} .$$
	
	Denote $y_1=-x_2$ and $y_2=-x_1$.  Then $-1 \leq y_2 \leq 1$, $y_1 > 1$, and 
	$$ \frac{ \left |  x_1 - x_2 - \sqrt{x_2^2-1}   \right | }{ \sqrt{x_1-x_2}} = \frac{ \left |  y_1 - y_2 - \sqrt{y_1^2-1}   \right | }{ \sqrt{y_1-y_2}} ,$$ 
	
	which is bounded by Case 4 above.

\vspace{.1in}
\noindent
Combining all the cases above, we have proved that $	\frac{ \left | Hf(x_1) - Hf(x_2)  \right | }{\left | x_1 - x_2 \right |^{\frac{1}{2}}} $ is bounded for all $x_1, x_2 \in \mathbb{R}$, so $Hf(x) \in C^{0,\frac{1}{2}}\left ( \mathbb{R} \right )$.

\vspace{.4in} 

\bibliographystyle{plain}
\bibliography{thesis}

\vspace{.2in} 

\end{document}